\documentclass
[11pt]{amsart}
\usepackage{amsmath,amsthm, amscd, amssymb, amsfonts}
\newtheorem{Lemma}{Lemma}
\newtheorem{Theorem}{Theorem}
\newtheorem{Prop}{Proposition}

\newtheorem{Cor}{Corollary}

\newtheorem{Ex}{Example}
\newcommand{\pf}{\medskip\noindent{\sc Proof: }}
\begin{document}
\title{Biproducts and Kashina's examples}
\author{David E. Radford}
\address{University of Illinois at Chicago \\
Department of Mathematics, Statistics and \\ Computer Science (m/c 240) \\
801 South Morgan Street \\
Chicago, IL   60608-7045} \email{radford@uic.edu}
\numberwithin{equation}{section}

\maketitle
\date{}
\begin{abstract}
{\small  \rm We revisit a class of examples described in the original paper on biproducts, expand the class, and provide a detailed analysis of the coalgebra and algebra structures of many of these examples. Connections with the semisimple Hopf algebras of dimension a power of two determined by Kashina are examined.

The finite-dimensional non-trivial semisimple cosemisimple Hopf algebras we construct are shown to be lower cosolvable. Some of these have one proper normal Hopf subalgebra and are not lower solvable.}
\end{abstract}
\setcounter{section}{0}
\section*{Introduction}\label{SecIntro}
Biproducts have a well-established position in the theory of Hopf algebras over a field $k$. They play a central role in the classification of pointed Hopf algebras \cite{AS2010} and arise quite often in the classification of semisimple Hopf algebras.
A finite-dimensional Hopf algebra is said to be trivial if it or its dual is a group algebra. There are two non-trivial semisimple Hopf algebras of dimension $12$ over an algebraically closed field whose characteristic is not $2$ or $3$; these are biproducts $B \times H$ with $H = k[\mathbb{Z}_4]$ or $H = k[\mathbb{Z}_2 \times \mathbb{Z}_2]$. In \cite[4.1]{RadProj} a method was given for constructing biproducts which are semisimple Hopf algebras. In this paper we expand on the method and describe the structure of the resulting Hopf algebras in great detail.

Let $B \times H$ be a biproduct. Then $H$ can be regarded as a Hopf subalgebra of $B \times H$ and such $B \times H$ is a free $H$-module with $\mathrm{Rank}_H(B \times H) = \mathrm{Dim}(B)$. Here we consider the rank $2$ case and determine all possibilities for $B$ with mild restrictions on $k$. Kashina's classification work \cite{KashYevg2000Dim16,KashYevg2003SS2fam2M,KashYevg2003SS2m,KashYevg2007} on semisimple Hopf algebras $A$ of dimension $2^n$ has produced many examples which have a group algebra $H' = k[G]$ of dimension $2^{n-1}$. For these examples we show that $A$ is a biproduct of the form $B \times H'$ only when $A \simeq k[\mathbb{Z}_2] \otimes H'$ as a Hopf algebra.

Throughout $k$ is a field and all vector spaces are over $k$. We denote the identity map of a vector space $V$ by $\mathrm{I}_V$, its dimension by $\mathrm{Dim}(V)$, and we denote the span of a subset $S$ of $V$ by $\mathrm{sp}(S)$. If $A$ is an algebra then $1$ or $1_A$ will denote its identity element. We denote the multiplicative group invertible elements (units) of $k$ by $k^\times$. Any one of \cite{Abe,DRS,LamRad,Mont,RBook,SweedlerBook} serves as a Hopf algebra reference for this paper. There is a very extensive literature on finite-dimensional semisimple Hopf algebras. The reader is directed to \cite{Masuoka1995SixEight,Masuoka1996Further,Mont1998,Natale1999PQSq,Natale2004PQR,Natale2007SemiSolv,Sommer2002} for starters.

A special feature of the biproducts in this paper. Throughout $H$ is a Hopf algebra over the field $k$ as are our biproducts $A = B \times H$. Generally the first factor $B$ has a $k$-algebra and a $k$-coalgebra structure which makes it a bialgebra only in a certain category. In this paper $B$, with these $k$-algebra and $k$-coalgebra structures, is a bialgebra over $k$ as well.
\section{Preliminaries}\label{SecPrelim}
For $r \geq 1$ let $\mathrm{C}_r(k)$ be the coalgebra over $k$ with basis $\{c_{i \, j}\}_{1 \leq i, j \leq r}$ whose structure is defined by
\begin{equation}\label{EqCrk}
\epsilon(c_{i \, j}) = \delta_{i, j} \;\; \mbox{and} \;\; \Delta (c_{i \, j}) = \sum_{\ell = 1}^r c_{i \, \ell} \otimes c_{\ell \, j}
\end{equation}
for all $1 \leq i, j \leq r$, where $\delta_{i, j}$ is the Kronecker delta function. The coalgebra $\mathrm{C}_r(k)$ is simple since $\mathrm{C}_r(k) \simeq \mathrm{M}_r(k)^*$, where $\mathrm{M}_r(k)$ is the algebra of $r \times r$ matrices over $k$. Let $C$ be any coalgebra over $k$. A subset $\{c_{i \, j}\}_{1 \leq i, j \leq r}$ of $C$ is said to satisfy the \textit{comatrix identities} if (\ref{EqCrk}) holds. We will find it convenient to use the indexing set $\{0, \ldots, r-1\}$ in place of $\{1, \ldots, r\}$ in several contexts. This is a harmless change of notation. When $(M, \rho)$ is a left $C$-comodule we write $\rho(m) = m_{(-1)} \otimes m_{(0)} \in C \otimes M$ for $m \in M$.

Now let $H$ be a bialgebra over $k$. A vector space $M$ is a left $H$-module under the trivial module action, meaning $h{\cdot}m = \epsilon(h)m$ for all $h \in H$ and $m \in M$. Likewise $M$ is a left $H$-comodule under the trivial comodule action, meaning $\rho(m) = 1 \otimes m$ for all $m \in M$.

A left $H$-module algebra is a left $H$-module $(B, \cdot)$, where $B$ is an algebra over $k$, such that $h{\cdot}1 =\epsilon(h)1$ and $h{\cdot}(bb') = (h_{(1)}{\cdot}b)(h_{(2)}{\cdot}b')$ for all $h \in H$ and $b, b' \in B$. When $(B, \cdot)$ is a left $H$-module algebra the tensor product of vector spaces $B \otimes H$ has an algebra structure, called the smash product, defined by $1_{B \otimes H} = 1_B \otimes 1_H$ and
\begin{equation}\label{EqSmash}
(b \otimes h)(b' \otimes h') = b(h_{(1)}{\cdot}b') \otimes h_{(2)}h'
\end{equation}
for all $b, b' \in B$ and $h, h' \in H$. The usual notation for this algebra is $B \# H$ and tensors $b \otimes h$ are written $b \# h$. Observe that $(b \# h)(b' \# h') = bb' \# hh'$ if either $h = 1$ or $b' = 1$. As a consequence the map $H \longrightarrow B \# H$ defined by $h \mapsto 1 \# h$ is a one-one algebra map. Identifying $H$ with a subalgebra of $B \# H$ in this way, we have that $B \# H$ is a free right $H$-module under multiplication with basis $\{b \# 1\}_{b \in S}$, where $S$ is a linear basis for $B$. Thus $\mathrm{Rank}_H(B \# H) = \mathrm{Dim}(B)$.

A left $H$-comodule algebra is a left $H$-comodule $(B, \rho)$, where $B$ is an algebra over $k$, such that $\rho(1) = 1_H \otimes 1$ and $\rho(bb') = b_{(-1)}b'_{(-1)} \otimes b_{(0)}b'_{(0)}$ for all $b, b' \in B$. A left $H$-module coalgebra is a left $H$-module $(C, \cdot)$, where $C$ is a coalgebra over $k$, such that $\epsilon (h{\cdot}c) = \epsilon(h)\epsilon(c)$ and $\Delta (h{\cdot}c) = h_{(1)}{\cdot}c_{(1)} \otimes h_{(2)}{\cdot}c_{(2)}$ for all $h \in H$ and $c \in C$.

A left $H$-comodule coalgebra is a left $H$-comodule $(C, \rho)$, where $C$ is a coalgebra over $k$, such that $c_{(-1)}\epsilon(c_{(0)}) = \epsilon(c)1$ and
$$
c_{(1)(-1)}c_{(2)(-1)} \otimes c_{(1)(0)} \otimes c_{(2)(0)} = c_{(-1)} \otimes c_{(0)(1)} \otimes c_{(0)(2)}
$$
for all $c \in C$. When $(C, \rho)$ is a left $H$-comodule coalgebra the tensor product of vector spaces $C \otimes H$ has a coalgebra structure, called the smash coproduct, defined by $\epsilon_{C \otimes H} = \epsilon_C \otimes \epsilon_H$ and $\Delta (c \otimes h) = (c_{(1)} \otimes c_{(2)(-1)}h_{(1)}) \otimes (c_{(2)(0)} \otimes h_{(2)})$ for $c \in C$. The usual notation for this coalgebra is $C \natural H$ and $c \otimes h$ is written $c \natural h$. The definitions of module, comodule algebra, coalgebra express what it means to be an algebra or coalgebra in the category of left $H$-modules or left $H$-comodules.

Biproducts are described in the context of a Yetter-Drinfel'd category ${}^H_H\mathcal{YD}$ which is defined for $H$. The (left-left) Yetter-Drinfel'd category, which is denoted by ${}^H_H\mathcal{YD}$, is described as follows. Objects are triples $(M, \cdot, \rho)$, where $(M, \cdot)$ is a left $H$-module and $(M, \rho)$ is a left $H$-comodule, such that
\begin{equation}\label{EqYDCondition}
h_{(1)}m_{(-1)} \otimes h_{(2)}{\cdot}m_{(0)} = (h_{(1)}{\cdot}m)_{(-1)}h_{(2)} \otimes (h_{(1)}{\cdot}m)_{(0)}
\end{equation}
for all $h \in H$ and $m \in M$. Morphisms are functions $f : M \longrightarrow N$ of objects which are left $H$-module and left $H$-comodule maps. Multiplication of morphisms is function composition.

Note that $k$ is an object of ${}^H_H\mathcal{YD}$ with the trivial module and comodule structures. Generally a vector space is an object of ${}^H_H\mathcal{YD}$ when equipped with the trivial left $H$-module and $H$-comodule structures. If $M, N$ are objects of ${}^H_H\mathcal{YD}$ then $M \otimes N$ is also where $h{\cdot}(m \otimes n) = h_{(1)}{\cdot}m \otimes h_{(2)}{\cdot}n$ and $\rho(m \otimes n) = m_{(-1)}n_{(-1)} \otimes  (m_{(0)} \otimes n_{(0)})$ for all $h \in H$, $m \in M$, and $n \in N$. For objects $M, N$ of ${}^H_H\mathcal{YD}$ there is a morphism $\sigma_{M, N} : M \otimes N \longrightarrow N \otimes M$ defined by $\sigma_{M, N}(m \otimes n) = m_{(-1)}{\cdot}n \otimes m_{(0)}$ for all $m \in M$ and $n \in N$. When $H$ is a Hopf algebra these morphisms are isomorphisms and can be thought of as taking the place of the ``twist" map of vector spaces.

An algebra of ${}^H_H\mathcal{YD}$ is a triple $(A, \cdot, \rho)$, where $A$ is a an algebra over $k$ whose structure maps are morphisms, that is $(A, \cdot)$ is a left $H$-module algebra and $(A, \rho)$ is a left $H$-comodule algebra. A coalgebra of ${}^H_H\mathcal{YD}$ is a triple $(C, \cdot, \rho)$, where $C$ is a $k$-coalgebra whose structure maps are morphisms, that is $(C, \cdot)$ is a left $H$-module coalgebra and $(C, \rho)$ is a left $H$-comodule coalgebra.

Let $A$ and $B$ be algebras of ${}^H_H\mathcal{YD}$. The tensor product of vector spaces $A \otimes B$ is an algebra of ${}^H_H\mathcal{YD}$, where $m_{A \otimes B} = (m_A \otimes m_B) {\circ} (\mathrm{I}_A \otimes \sigma_{B, A}\otimes \mathrm{I}_B)$ defines the product of $A \otimes B$ in terms of those of $A$ and $B$. The convention is to write $a \underline{\otimes} b$ for a tensor $a \otimes b$ and denote this algebra by $A \underline{\otimes} B$. Thus
$$
(a \underline{\otimes} b)(a' \underline{\otimes} b') = a(b_{(-1)}{\cdot}a') \underline{\otimes} b_{(0)}b'
$$
for all $a, a' \in A$ and $b, b' \in B$. The unit of $A \underline{\otimes} B$ is given by $1_{A \underline{\otimes} B} = 1_A \underline{\otimes}1_B$.

Let $C$ and $D$ be coalgebras of ${}^H_H\mathcal{YD}$. The tensor product of vector spaces $C \otimes D$ is a coalgebra of ${}^H_H\mathcal{YD}$,  where $\Delta_{C \otimes D} = (\mathrm{I}_C \otimes \sigma_{C, D}\otimes \mathrm{I}_D) {\circ} (\Delta_C \otimes \Delta_D)$.  The convention is to write $c \underline{\otimes} d$ for a tensor $c \otimes d$ and denote this coalgebra by $C \underline{\otimes} D$. Thus
$$
\Delta (c \underline{\otimes} d) = (c_{(1)} \underline{\otimes} c_{(2)(-1)}{\cdot}d_{(1)}) \otimes (c_{(2)(0)} \underline{\otimes} d_{(2)})
$$
for all $c \in C$ and $d \in D$. The counit is given by $\epsilon_{C \underline{\otimes} D} = \epsilon_C \otimes \epsilon_D$. A bialgebra $A$ of ${}^H_H\mathcal{YD}$ is an algebra and coalgebra of ${}^H_H\mathcal{YD}$ whose coalgebra structure maps are algebra morphisms, which is to say $\epsilon$ and $\Delta : A \longrightarrow A \underline{\otimes} A$ are morphisms of ${}_H^H\mathcal{YD}$ and $k$-algebra maps; in particular
\begin{equation}\label{EqAtenAUnderDelta}
\Delta (bb') = b_{(1)}(b_{(2)(-1)}{\cdot}b'_{(1)}) \otimes b_{(2)(0)}b'_{(2)}
\end{equation}
for all $b, b'  \in B$. A Hopf algebra $B$ of ${}^H_H\mathcal{YD}$ bialgebra $B$ of ${}^H_H\mathcal{YD}$ with a morphism $S : B \longrightarrow B$ such that $S(b_{(1)})b_{(2)} = \epsilon(b)1 = b_{(1)}S(b_{(2)})$ for all $b \in B$. We note that any linear endomorphism of $B$ which satisfies the two preceding equations is a morphism of ${}^H_H\mathcal{YD}$.

Suppose that $H$ is a Hopf algebra over $k$ and $B$ is a Hopf algebra of ${}^H_H\mathcal{YD}$. Then the vector space $A = B \otimes H$ is a Hopf algebra over $k$, called the biproduct of $B$ and $H$, with the smash product algebra and smash coproduct coalgebra structures.

We denote the biproduct $A$ by $B \times H$. By virtue of our discussion of smash products, $B \times H$ is a free right $H$-module and $\mathrm{Rank}_H(B \times H) = \mathrm{Dim}(B)$. Various other notations for tensors $a \otimes h$ of $A$ are sometimes used, namely $a \times h$, $a \# h$, or $a \natural h$, when a structure is to be emphasized: biproduct, smash product, or smash coproduct.

We note a biproduct $B \times H$ has the tensor product algebra structure if and only if the $H$-module action is trivial, and that the biproduct is commutative if and only if the $H$-module action is trivial and $B$, $H$ are commutative. Likewise, $B \times H$ has the tensor product coalgebra structure if and only if the $H$-comodule action is trivial, and the biproduct is cocommutative if and only if the $H$-comodule action is trivial and $B$, $H$ are cocommutative.

Suppose $A = B \times H$ and $H$ are Hopf algebras over $k$. Then there is a Hopf algebra projection $A$ onto $H$, meaning there are Hopf algebra maps $\pi : A \longrightarrow H$ and $j : H \longrightarrow A$ such that $\pi {\circ} j = \mathrm{I}_H$. Conversely, if $A$ and $H$ are Hopf algebras over $k$ and there are such maps then $A \simeq B \times H$ for some Hopf algebra $B$ of ${}^H_H\mathcal{YD}$.

A minor but useful observation. Algebras, coalgebras, bialgebras, and Hopf algebras over $k$ are algebras, coalgebras, bialgebras, and Hopf algebras respectively of ${}^H_H\mathcal{YD}$ when they are equipped with the trivial left $H$-module and $H$-comodule structures.

Biproducts $B \times H$ were defined and studied in \cite{RadProj} in connection with Hopf algebras with a projection. In \cite{Majid1993Skly,Majid1994Braided} it was realized that $B$ is a Hopf algebra in the braided monoidal category ${}_H^H\mathcal{YD}$ and that this is the proper theoretical setting for $B$. The category ${}_H^H\mathcal{YD}$ is the category $A$-$\mathbf{cbm}$ defined and studied in \cite{Yetter1990}, where $H = A$ and objects are referred to as left crossed bimodules over $A$. Terminology for ${}_H^H\mathcal{YD}$, and study of variations of this category, is found in \cite{Rad-Tow1993}.

We review the left, right adjoint module actions $\mathrm{ad}_\ell, \mathrm{ad}_r : H \otimes H \longrightarrow H$ of a Hopf algebra $H$ with antipode $S$ on itself for the last sections. By definition $h{\succ}a = \mathrm{ad}_\ell(h \otimes a) = h_{(1)}aS(h_{(2)})$ and $a {\prec}h = \mathrm{ad}_r(a \otimes h) = S(h_{(1)})ah_{(2)}$ for all $h, a \in H$. Observe that $ha = (h_{(1)}{\succ}a)h_{(2)}$ and $ah = h_{(1)}(a{\prec}h_{(2)})$ for all $h, a \in H$. Thus is $M$ is a submodule of $H$ under the left and right adjoint actions $HM = MH$.

A Hopf subalgebra $K$ of $H$ is said to normal if it is a submodule under left and right adjoint actions. In this case $K^+ = \mathrm{Ker}(\epsilon) \cap K$ is a submodule under the actions and therefore $HK^+ = K^+H$. Thus the (Hopf) ideal of $H$ generated by $K^+$ is the left ideal, or the right ideal, of $H$ generated by $K^+$.

\section{Projections onto Hopf algebras of rank $2$.}\label{SecRank2}
Throughout this section the field $k$ has characteristic not $2$ and $H$ is a Hopf algebra over $k$. In this section we determine all $2$-dimensional Hopf algebras $B$ in ${}_H^H\mathcal{YD}$. There is one and only one which is semisimple, namely the group algebra $k[\mathbb{Z}_2]$ over $k$ with trivial $H$-module and $H$-comodule actions.

Let $B$ be a $2$-dimensional left $H$-module with basis $\{1, n\}$. Let $\alpha, \beta \in H^*$ be determined by
\begin{equation}\label{EqhnBetaAlpha}
h{\cdot}n = \beta(h)1 + \alpha(h)n
\end{equation}
for all $h \in H$. Suppose that
\begin{equation}\label{EqH1Epsilon1}
h{\cdot}1 = \epsilon(h)1
\end{equation}
for all $h \in H$. The preceding equation holds if $B$ is an $H$-module algebra and $1 = 1_B$. Since $B$ is a left $H$-module
\begin{equation}\label{EqAlpha}
\alpha \in \mathrm{Alg}_k(H, k) = \mathrm{G}(H^o)
\end{equation}
and $\beta(hh') = \epsilon(h)\beta(h') + \beta(h)\alpha(h')$ for all $h, h' \in H$. Since $H^o$ is a Hopf algebra $\alpha$ is invertible.

It will be convenient to place restrictions on $n$ when $B$ is an algebra. Suppose $B$ is an algebra over $k$ with basis $\{1, m\}$. Then $m^2 + rm + s1 = 0$ for some $r, s \in k$. Since the characteristic of $k$ is not $2$ there is a basis $\{1, n\}$ for $B$ such that $n^2 = \varpi 1$ for some $\varpi \in k$, namely with $n = m + (r/2)1$.

Now assume $B$ is a $2$-dimensional left $H$-module algebra. Then $h{\cdot}1 = \epsilon(h)1$ for all $h \in H$. Choose a basis $\{1, n\}$ for $B$ such that $n^2 = \varpi1$, where $\varpi \in k$, and let $\alpha, \beta \in H^*$ be determined by (\ref{EqhnBetaAlpha}). Using the module algebra axiom $h{\cdot}(bb') = (h_{(1)}{\cdot} b)(h_{(2)}{\cdot} b')$ for all $h \in H$ and $b, b' \in B$ we compute
\begin{eqnarray*}
h{\cdot}n^2 & = & (h_{(1)}{\cdot}n)(h_{(2)}{\cdot}n) \\
& = & (\beta(h_{(1)})1 + \alpha(h_{(1)})n)(\beta(h_{(2)})1 + \alpha(h_{(2)})n) \\
& = & \beta^2(h)1 + (\beta\alpha + \alpha\beta)(h)n + \alpha^2(h)n^2 \\
& = & (\beta^2 + \varpi\alpha^2)(h)1 + (\beta\alpha + \alpha\beta)(h)n.
\end{eqnarray*}
We have shown:
\begin{Prop}\label{Prop1n}
Let $H$ be a Hopf algebra over the field $k$ and let $B$ be a left $H$-module algebra with basis $\{1, n\}$ such that $n^2 = \varpi1$, where $\varpi \in k$. Let $\alpha, \beta \in H^*$ be determined by (\ref{EqhnBetaAlpha}). Then:
\begin{enumerate}
\item[{\rm (a)}]
$\beta\alpha = - \alpha\beta$ \, \mbox{and}
\item[{\rm (b)}]
$\beta^2 = \varpi(\epsilon - \alpha^2)$.
\end{enumerate}
\qed
\end{Prop}

Now suppose that $B$ is a left $H$-module coalgebra, where $B$ has basis $\{1, n\}$, and whose coproduct satisfies
\begin{equation}\label{EqGB}
1 \in \mathrm{G}(B) \;\; \mbox{and} \;\; h{\cdot}1 = \epsilon(h)1
\end{equation}
for all $h \in H$. This is the case when $B$ is a bialgebra of ${}^H_H\mathcal{YD}$ and $1 = 1_B$. Write
$$
\Delta(n) = r \otimes 1 + s \otimes n,
$$
where $r, s \in B$. The coalgebra axioms applied to $n$ yield
$$
\Delta(r) = r \otimes 1 + s \otimes r, \;\;  \epsilon(r) = 0, \;\; \mbox{and}
$$
\begin{equation}\label{EqDeltas}
s \in \mathrm{G}(B).
\end{equation}
From the first coproduct formula we deduce that $n = r + \epsilon(n)s$
and therefore
\begin{equation}\label{EqDeltans1}
\Delta(n) = n \otimes 1 + s \otimes n - \epsilon(n)s \otimes 1
\end{equation}
follows from the second. Since $\{1, s\} \subseteq \mathrm{G}(B)$ either $\{1, s\}$ is a basis for $B$ or $s = 1$. Using the module coalgebra axiom $h{\cdot}\Delta(n) = \Delta(h{\cdot}n)$ for all $h \in H$ and (\ref{EqDeltans1}) two equations follow:
$$
(\beta \rightharpoonup h){\cdot}s  - \epsilon(n)h{\cdot}s = -\epsilon(n)\alpha(h)s
$$
and
$$
(\alpha \rightharpoonup h){\cdot}s = \alpha(h)s
$$
for all $h \in H$. Since $\alpha$ is invertible, the second equation is equivalent to
\begin{equation}\label{Eqhs}
h{\cdot}s = \epsilon(h)s
\end{equation}
for all $h \in H$. Since $\epsilon(h{\cdot}b) = \epsilon(h)\epsilon(b)$ for all $h \in H$ and $b \in B$, applying $\epsilon$ to both sides of the first gives
$\beta = \epsilon(n)(\epsilon - \alpha)$. In particular $\beta$ and $\alpha$ commute.
\begin{Cor}
Assume the hypothesis of Proposition \ref{Prop1n} and also that $B$ is a left $H$-module coalgebra such that (\ref{EqGB}) holds. Then $\beta = 0$; therefore
\begin{equation}\label{EqHDotN}
h{\cdot}n = \alpha(h)n
\end{equation}
for all $h \in H$.
\end{Cor}

\pf
$\beta\alpha = -\alpha\beta$ by part (a) of Proposition \ref{Prop1n}. Thus $\beta = 0$ since $\alpha$ and $\beta$ commute, $\alpha$ is invertible, and the characteristic of $k$ is not $2$. The last assertion now follows by (\ref{EqhnBetaAlpha}).
\qed
\medskip

Now suppose $B$ is a left $H$-comodule coalgebra with basis $\{1, n\}$ where
\begin{equation}\label{EqRho1}
\rho(1) = 1_H \otimes 1 \;\;\mbox{and} \;\; 1 \in \mathrm{G}(B),
\end{equation}
We observe that (\ref{EqRho1}) is satisfied when $B$ is a bialgebra of ${}_H^H\mathcal{YD}$ and $1 = 1_B$. Write
$$
\rho(n) = x \otimes 1 + y \otimes n,
$$
where $x, y \in B$. The comodule axioms yield
\begin{equation}\label{EqyGH}
y \in \mathrm{G}(H),
\end{equation}
$$
\Delta (x) = x \otimes 1 + y \otimes x, \;\; \mbox{and} \;\; \epsilon(x) = 0.
$$
The comodule coalgebra axiom $b_{(-1)}\epsilon(b_{(0)}) =\epsilon(b)1$ for all $b \in B$ applied to $b = n$ gives $x = \epsilon(n) (1 - y)$; in particular $x$ and $y$ commute. Suppose further that the coalgebra structure of $B$ satisfies (\ref{EqDeltas}) and (\ref{EqDeltans1}). The comodule coalgebra axiom
$$
b_{(1)(-1)}b_{(2)(-1)} \otimes b_{(1)(0)} \otimes b_{(2)(0)} = b_{(-1)} \otimes b_{(0)(1)} \otimes b_{(0)(2)}
$$
for all $b \in B$ applied to $b = n$ yields
$
s_{(-1)}y \otimes s_{(0)} = y \otimes s.
$
Now $y$ is invertible since $y \in \mathrm{G}(H)$ and $H$ is a Hopf algebra. Thus the preceding equation implies
\begin{equation}\label{EqRhos}
\rho(s) = 1 \otimes s.
\end{equation}

Now suppose that $B$ is also a left $H$-comodule algebra and further assume $1 = 1_B$ and $n^2 = \varpi1$, where $\varpi \in k$. Applying the comodule algebra axiom
$$
\rho(ab) = a_{(-1)}b_{(-1)} \otimes a_{(0)}b_{(0)}
$$
to $a = b = n$ we obtain
$$
\rho (n^2) = (x^2 + \varpi y^2) \otimes 1 + (xy + yx) \otimes n.
$$
Therefore $xy = -yx$ and $x^2 =  \varpi (1 - y^2)$. We have noted that $xy = yx$ and $y$ is invertible. Thus $x = 0$; in particular
\begin{equation}\label{EqRhoNY}
\rho(n) = y \otimes n.
\end{equation}

We are close to determining all $2$-dimensional Hopf algebras in ${}_H^H\mathcal{YD}$. Suppose that $B$ is a $2$-dimensional algebra and coalgebra over $k$ and has a left $H$-module and left $H$-comodule structure which affords it the structure of a left $H$-module algebra and coalgebra and left $H$-comodule algebra and coalgebra. Then there is a basis $\{1, n\}$ for $B$ such that $n^2 = \varpi 1$ for some $\varpi \in k$. Choose such a basis. Then (\ref{EqAlpha})--(\ref{EqRhoNY}) hold. Next we establish a necessary and sufficient condition for $B \in {}_H^H\mathcal{YD}$, that is for (\ref{EqYDCondition}) to hold for all $h \in H$ and $b \in B$. It is easy to see that (\ref{EqYDCondition}) holds for $m = 1$ and all $h \in H$. For $m = n$ observe that (\ref{EqYDCondition}) holds for all $h \in H$ if and only if $(\alpha \rightharpoonup h)y = y(h \leftharpoonup \alpha)$, or equivalently
\begin{equation}\label{EqAlphahy}
\alpha \rightharpoonup h \leftharpoonup \alpha^{-1} = yhy^{-1},
\end{equation}
for all $h \in H$.
Assume (\ref{EqAlphahy}) holds as well and that $B$ is a bialgebra in ${}_H^H\mathcal{YD}$. There are two cases to consider.
\medskip

\noindent
\underline{Case} 1: $s \neq 1$.

Here $\{1, s\}$ is a basis for $B$. By (\ref{EqGB}) and (\ref{Eqhs}) we have $h{\cdot}b = \epsilon(h)b$ for all $h \in H$ and $b \in B$. By (\ref{EqRho1}) and (\ref{EqRhos}) we have $\rho(b) = 1 \otimes b$ for all $b \in B$. Therefore the $H$-module and comodule actions are trivial; in particular $\alpha = \epsilon$ and $y = 1$. Since the actions are trivial $B$ is a bialgebra over $k$ and is therefore the group algebra $k[\mathbb{Z}_2]$ or the semigroup algebra $k[S]$, where $S = \{1, s \}$ is the semigroup with neutral element $1$ and $s^2 = s$.
\medskip

\noindent
\underline{Case} 2: $s = 1$.

Then (\ref{EqAtenAUnderDelta}) for $b = b' = n$ is
$$
\Delta (n^2) = n^2 \underline{\otimes} 1 + 1 \underline{\otimes} n^2 + (1 + \alpha(y))n \underline{\otimes} n - 2\epsilon(n)(n \underline{\otimes}1 + 1 \underline{\otimes} n) - \epsilon(n)^2 1 \underline{\otimes} 1.
$$
Since $n^2 = \varpi 1$ and the characteristic of $k$ is not $2$ we conclude
$$
1 + \alpha(y) = 0 = \epsilon (n)
$$
The last equation implies $\varpi = 0$, hence $n^2 = 0$.

We encompass Case 2 in the following structure. For $\alpha \in \mathrm{G}(H^o)$ and $y \in \mathrm{G}(H)$ let $B_{\alpha, y}$ be the $k$-algebra with basis $\{1, n\}$ where $n^2 = 0$. Give  $B = B_{\alpha, y}$ the $k$-coalgebra structure determined by
$$
\Delta(1) = 1 \otimes 1 \;\; \mbox{and} \;\; \Delta (n) = n\otimes 1 + 1 \otimes n,
$$
the left $H$-module structure determined by
$$
h{\cdot}1 = \epsilon(h)1 \;\; \mbox{and} \;\; h{\cdot}n = \alpha(h)n
$$
for all $h \in H$, and the left $H$-comodule structure determined by
$$
\rho(1) = 1 \otimes 1 \;\; \mbox{and} \;\; \rho(n) = y \otimes n.
$$
We leave the reader with the easy exercise of showing that $B_{\alpha, y}$ is a bialgebra of ${}_H^H\mathcal{YD}$ if and only if $\alpha(h) = -1$ and (\ref{EqAlphahy}) holds for all $h \in H$. For any $\alpha, y$ observe that the endomorphism $S$ of $B_{\alpha, y}$ determined by $S(1) = 1$ and $S(n) = -n$ is an antipode for $B_{\alpha, y}$. Collecting results:
\begin{Theorem}\label{TheoremMain}
Let $H$ be a Hopf algebra over the field $k$ and suppose that the characteristic of $k$ is not $2$. Then:
\begin{enumerate}
\item[{\rm (a)}] Suppose $\alpha \in \mathrm{G}(H^o)$ and $y \in \mathrm{G}(H)$ satisfy $\alpha(y) = -1$ and (\ref{EqAlphahy}) for all $h \in H$. Then $B_{\alpha, y}$ is a Hopf algebra in ${}_H^H\mathcal{YD}$.
\item[{\rm (b)}] Let $B$ be a $2$-dimensional Hopf algebra in ${}_H^H\mathcal{YD}$. Then $B \simeq k[\mathbb{Z}_2]$ with its $k$-Hopf algebra structure and trivial $H$-module and comodule structures, or $B \simeq B_{\alpha, y}$, where $\alpha \in \mathrm{G}(H^o)$ and $y \in \mathrm{G}(H)$ satisfy $\alpha(y) = -1$ and (\ref{EqAlphahy}) for all $h \in H$.
\end{enumerate}
\qed
\end{Theorem}

\begin{Cor}\label{Cor1}
Let $A$ be a finite-dimensional Hopf algebra over the field $k$ and suppose that the characteristic of $k$ is not $2$. Assume that $H$ is Hopf subalgebra of $A$, there is a Hopf algebra projection $A$ onto $H$, and that $A$ is a free right $H$-module of rank $2$. Then:
\begin{enumerate}
\item[{\rm (a)}] $A \simeq B_{\alpha, y} \times H$ for some $\alpha \in \mathrm{G}(H^o)$ and $y \in \mathrm{G}(H)$ which satisfy $\alpha(y) = -1$ and (\ref{EqAlphahy}) for all $h \in H$, or
\item[{\rm (b)}] $A \simeq k[\mathbb{Z}_2] \otimes H$ as Hopf algebras over $k$. The latter is the case when $A$ is semisimple.
\end{enumerate}
\end{Cor}

\pf
If part (a) does not hold then $A = B \otimes H$ as a bialgebra. Therefore grouplike elements of $B$ must be invertible which means $B \simeq k[\mathbb{Z}_2]$. As for the semisimple assertion, by part (b) of \cite[\S 2 Proposition 3]{RadProj} if a biproduct $B \times H$ is a semisimple Hopf algebra there is a $\Lambda \in $B such that $b\Lambda = \epsilon(b)\Lambda$ for all $b \in B$ and $\epsilon(\Lambda) = 1$. There is no such element in $B_{\alpha, y}$.
\qed
\medskip

Many of the non-trivial semisimple Hopf algebras $A$ which Kashina has described \cite{KashYevg2000Dim16,KashYevg2007} have a Hopf subalgebra $H$ which is a group algebra and $A$ is a free right $H$-module of rank $2$. None of them has a Hopf algebra projection onto $H$, for:
\begin{Cor}\label{Cor2}
Under the hypothesis of the preceding corollary, when $A$ is semisimple:
\begin{enumerate}
\item[{\rm (a)}] If $H$ is a group algebra then $A$ is a group algebra.
\item[{\rm (b)}] If $H^*$ is a group algebra then $A^*$ is a group algebra.
\item[{\rm (c)}] If $H$ is trivial then $A$ is trivial.
\end{enumerate}
\end{Cor}

\pf
Since the characteristic of $k$ is not $2$ we have $(k[\mathbb{Z}_2] \otimes H)^* \simeq k[\mathbb{Z}_2]^* \otimes H^* \simeq k[\mathbb{Z}_2] \otimes H^*$ as Hopf algebras. Thus Corollary \ref{Cor2} follows by the preceding corollary.
\qed
\section{A basic construction revisited}\label{SecConstRev}
In this section we reconsider the construction of \cite[4.2]{RadProj}. We use our results to produce several families of non-trivial semisimple cosemisimple Hopf algebras $A$ with a projection onto a Hopf subalgebra $H$. As an example, for $n \geq 3$ and odd, and $H = k[{\mathbb Z}_2 \times {\mathbb Z }_2]$, we construct such a biproduct $A = B \times H$ where $A$ is a free right $H$-module of rank $n$.

We loosen the assumptions of \cite[4.2]{RadProj} a bit. Let $G$ and $\mathcal{G}$ be any groups and $B = k[\mathcal{G}]$. Regard group algebras as Hopf algebras in the usual way. We make the identification $\mathrm{Aut}_\mathrm{Bialg}(B) = \mathrm{Aut}_\mathrm{Group}(\mathcal{G})$ via restriction to $\mathcal{G}$.

Suppose $\pi : G \longrightarrow \mathrm{Aut}_\mathrm{Group}(\mathcal{G})$ is a group homomorphism and set $H = k[G]$. Then $B$ is a left $H$-module where $g{\cdot}b = \pi(g)(b)$ for all $g \in G$ and $b \in \mathcal{G}$. By parts (a)(i) and (b)(i) of \cite[\S 4 Lemma 1]{RadProj} we conclude that $B$ is a left $H$-module algebra and a left $H$-module coalgebra.

Suppose further that $\theta \in \mathrm{Aut}_\mathrm{Group}(\mathcal{G})$ and that the $(\theta)$-orbits of $\mathcal{G}$ are finite. Regard $\theta$ as a bialgebra automorphism of $B$. Then the set $S$ of eigenvalues of $\theta$ in $k$ consists of roots of unity. Furthermore $\theta$ is a split operator, meaning
\begin{equation}\label{EqSplitOerator}
B = \bigoplus_{\lambda \in k}B_{(\lambda)} \;\; ( = \bigoplus_{\lambda \in S}B_{(\lambda)}),
\end{equation}
where $B_{(\lambda)}$ is the subspace of all $b \in B$ which satisfy $(\theta - \lambda \mathrm{I}_B)^n(b) = 0$ for some $n \geq 0$. By parts (a) and (b) of \cite[4.1]{RadProj} we have
\begin{equation}\label{EqSp1}
1 \in B_{(1)},
\end{equation}
\begin{equation}\label{EqSp2}
B_{(\lambda)}B_{(\lambda')} \subseteq B_{(\lambda \lambda')} \;\; \mbox{for all $\lambda, \lambda' \in k$, and}
\end{equation}
\begin{equation}\label{EqSp3}
\epsilon (B_{(\lambda)}) = (0) \;\; \mbox{unless $\lambda = 1$ and}
\end{equation}
\begin{equation}\label{EqSp4}
\Delta(B_{(\lambda)}) \subseteq \sum_{\lambda'\lambda'' = \lambda} B_{(\lambda')} \otimes B_{(\lambda'')} \;\; \mbox{for all $\lambda \in k$.}
\end{equation}

Now let $U$ be the subgroup of $k^\times$ generated by $S$. Assume that $U$ is identified with a subgroup $\mathbf{U}$ of $G$ which satisfies
$$
\mathbf{U} \subseteq \mathrm{Ker}(\pi) \cap \mathrm{Z}(G)
$$
and denote this identification of $U$ and $\mathbf{U}$ by $\lambda \mapsto$ \mbox{ \boldmath $\lambda$}.

For $g \in G$ we set $B_g = (0)$ unless $g =$ \mbox{\boldmath $\lambda$} for some \mbox{\boldmath $\lambda$} $\in \mathbf{U}$ in which case we set $B_g = B_{(\lambda)}$. Then
\begin{equation}\label{EqBgSum}
B = \bigoplus_{g \in G}B_g
\end{equation}
by (\ref{EqSplitOerator}).

Now assume that
\begin{equation}\label{EqCommute}
\pi(g) \circ \theta = \theta \circ \pi(g)
\end{equation}
for all $g \in G$. Then $\pi (g)(B_{(\lambda)}) \subseteq B_{(\lambda)}$ for all $g \in G$ and $\lambda \in k$. Therefore $B_g$ is a left $H$-submodule of $B$ for all $g \in G$. Since $B_g = (0)$ unless $g \in \mathrm{Z}(G)$, $h{\cdot}B_g \subseteq B_g = B_{hgh^{-1}}$ for all $h, g \in G$. Regard $B$ as a left $H$-comodule according to  $\rho (b) = g \otimes b$ for all $g \in G$ and $b \in B_g$. By parts (a)(ii) and (b)(ii) of \cite[\S 4 Lemma 1]{RadProj} we conclude that $B$ is a left $H$-comodule algebra and a left $H$-comodule coalgebra. By part (c) of the same $B$ with its $H$-module and comodule structures is an object of ${}^H_H\mathcal{YD}$. Since $B_g \neq (0)$ implies $g \in \mathrm{Ker}(\pi)$, by \cite[\S 4, Proposition 5]{RadProj} it follows that $B$ is a bialgebra of ${}_H^H\mathcal{YD}$, hence is a Hopf algebra of the same. Therefore the biproduct $A = B \times H$ is defined.

\textit{For the remainder of this section we assume that whenever} $r$ \textit{is the length of a} $(\theta)$-\textit{orbit of} $\mathcal{G}$ \textit{then} $k$ \textit{contains a primitive} $r^{th}$ \textit{root of unity}. We will show that $A$ is cosemisimple by showing it is the direct sum of comatrix coalgebras and we will describe them explicitly.

Let $b \in \mathcal{G}$. By assumption the $(\theta)$-orbit $\mathcal{O}_b$ of $b$ is finite. Let $r = |\mathcal{O}_b|$. Then $\mathcal{O}_b = \{b, \theta(b), \ldots, \theta^{r-1}(b)\}$ and $\theta^r(b) = b$. Let $\omega \in k$ be an $r^{th}$ root of unity and set
$$
b_\omega = \sum_{i = 0}^{r-1} (\omega^{-i}/r)\theta^i(b).
$$
Then $\theta(b_{\omega}) = \omega b_\omega$. By assumption $k$ contains a primitive $r^{th}$ root of unity $\lambda$. Therefore the span of $\mathcal{O}_b$ has basis of eigenvectors $\{b_{\lambda^0}, b_{\lambda^1}, \ldots, b_{\lambda^{r-1}}  \}$. We next compute the $\theta^i(b)$'s in terms of this basis.

Fix $0 \leq i \leq r-1$. Then $\theta^i(b) = \sum_{\ell = 0}^{r-1} \alpha_{i \, \ell} b_{\lambda^\ell}$ for some $\alpha_{i \, \ell} \in k$. This equation can be expressed as
$$
\theta^i(b) = \sum_{\ell, j = 0}^{r-1} \alpha_{i \, \ell}(\lambda^{-\ell j}/r)\theta^j(b)
$$
which is equivalent to
$$
\sum_{\ell = 0}^{r-1} \alpha_{i \, \ell}(\lambda^{- \ell j}/r) = \delta_{i, j}
$$
for all $0 \leq i, j \leq r-1$. Let $\mathcal{A}, \mathcal{V} \in \mathrm{M}_r(k)$ be defined by $\mathcal{A} = (\alpha_{i \, j})$ and $\mathcal{V} = (\nu_{i \, j})$, where $\nu_{i \, j} = \lambda^{-i j}/r$. The previous set of equations may be expressed $\mathcal{A}\mathcal{V} = \mathrm{I}_r$. Thus $\mathcal{A}$ and $\mathcal{V}$ are inverses. At this point it is easy to see that $\alpha_{i \, j} = \lambda^{i j}$. We have shown
\begin{equation}\label{EqThetaIb}
\theta^i(b) = \sum_{\ell = 0}^{r-1} \lambda^{i \ell} b_{\lambda^\ell}
\end{equation}
for all $0 \leq i \leq r-1$.

Let $\mathrm{C}(\mathcal{O}_b, \mbox{\boldmath $\lambda$})$ be the span of the $\theta^i(b) \times$ \mbox{\boldmath $\lambda$}${}^j$'s, where $0 \leq i, j \leq r-1$. Then $\mathrm{Dim}(\mathrm{C}(\mathcal{O}_b, \mbox{\boldmath $\lambda$})) = r^2$. We will show that $\mathrm{C}(\mathcal{O}_b, \mbox{\boldmath $\lambda$})$ is a subcoalgebra of $A$ which is isomorphic to the comatrix coalgebra $\mathrm{C}_r(k)$.

Let $v_i = \theta^i(b) \times 1$ for $0 \leq i \leq r-1$. We use the fact that $\theta^i(b)$ is grouplike and (\ref{EqThetaIb}) to calculate
\begin{eqnarray*}
\Delta (v_i) & = & (\theta^i(b) \times \theta^i(b)_{(-1)}) \otimes (\theta^i(b)_{(0)} \times 1) \\
& = & \sum_{\ell = 0}^{r-1} (\theta^i(b) \times \lambda^{i \ell}b_{\lambda^\ell \, (-1)}) \otimes (b_{\lambda^\ell \, (0)} \times 1) \\
& = & \sum_{\ell, j = 0}^{r-1} (\theta^i(b) \times \lambda^{i \ell}\mbox{\boldmath{$\lambda$}${}^\ell$}) \otimes ((\lambda^{-\ell j}/r) \theta^j(b) \times 1) \\
& = & \sum_{j = 0}^{r-1} c_{i \, j} \otimes v_j,
\end{eqnarray*}
where
$$
c_{i \, j} = \sum_{\ell = 0}^{r-1} \theta^i(b) \times (\lambda^{\ell (i - j)}/r) \mbox{\boldmath{$\lambda$}${}^\ell$}.
$$
Since $\{v_0, \ldots, v_{r-1} \}$ is linearly independent the $c_{i \, j}$'s satisfy the comatrix identities. Note that the $c_{i \, j}$'s belong to $\mathrm{C}(\mathcal{O}_b, \mbox{\boldmath $\lambda$})$. To justify our assertions about $\mathrm{C}(\mathcal{O}_b, \mbox{\boldmath $\lambda$})$ we need only show that the $c_{i \, j}$'s form a linearly independent set.

Suppose $\sum_{i, j = 0}^{r-1} x_{i \, j} c_{i \, j} = 0$, where $x_{i \, j} \in k$. Then
$$
\sum_{i, j, \ell = 0}^{r-1} x_{i \, j}\lambda^{\ell (i - j)} \theta^i(b) \times \mbox{\boldmath $\lambda$}^\ell = 0,
$$
or equivalently
$$
0 = \sum_{j = 0}^{r-1} x_{i \, j}\lambda^{\ell (i - j)} = \lambda^{i \ell}(\sum_{j = 0}^{r-1} x_{i \, j}\lambda^{-j \ell}),
$$
and thus
$$
\sum_{j = 0}^{r-1} x_{i \, j}\lambda^{-j \ell} = 0,
$$
for all $0 \leq i, \ell \leq r-1$. Let $\mathcal{X} = (x_{i \, j})$. The last set of equations can be expressed as $\mathcal{X}(r\mathcal{V}) = 0$ from which $\mathcal{X} = 0$ follows since $r\mathcal{V}$ is invertible. We have shown that the $c_{i \, j}$'s form a linearly independent set. Thus $\mathrm{C}(\mathcal{O}_b, \mbox{\boldmath $\lambda$})$ is a subcoalgebra of $A$ and $\mathrm{C}(\mathcal{O}_b, \mbox{\boldmath $\lambda$}) \simeq \mathrm{C}_r(k)$.

Let $g \in G$. Then $1 \times g \in \mathrm{G}(A)$ and $(c \times h)(1 \times g) = c \times hg$ for all $c \in B$ and $h \in H$. Therefore the coalgebra $\mathrm{C}(\mathcal{O}_b, \mbox{\boldmath $\lambda$})(1 \times g) \simeq \mathrm{C}_r(k)$ also and has basis $\{\theta^i(b) \times \mbox{\boldmath{$\lambda$}}{}^jg \}_{0 \leq i, j \leq r-1}$.

Let $\mathbf{U}_r$ be the cyclic subgroup of $G$ generated by $\mbox{\boldmath $\lambda$}$. Then $\mathbf{U}_r$ is the only subgroup of $\mathbf{U}$ of order $r$. Let $\mathcal{O} = \mathcal{O}_b$. Then $\mathcal{O} \times \mathbf{U}_r$ is a basis for $\mathrm{C}(\mathcal{O}_b, \mbox{\boldmath $\lambda$})$. We set $\mathrm{C}(\mathcal{O}, \mathbf{U}_r) = \mathrm{C}(\mathcal{O}_b, \mbox{\boldmath $\lambda$})$. It will be useful to note
\begin{equation}\label{RhoB}
\rho(\mathrm{sp}(\mathcal{O}_b)) \subseteq k[\mathbf{U}_r] \otimes \mathrm{sp}(\mathcal{O}_b),
\end{equation}
where $r = |\mathcal{O}_b|$.

Now let $R(G, \mathbf{U}_r)$ be a set of right coset representatives of $\mathbf{U}_r$ in $G$. Then
\begin{equation}\label{EqSS1}
\mathrm{sp}(\mathcal{O}) \times H = \bigoplus_{g \in R(G, \mathbf{U}_{|\mathcal{O}|})}  \mathrm{C}(\mathcal{O}, \mathbf{U}_{|\mathcal{O}|})(1 \times g) \simeq \bigoplus_{g \in R(G, \mathbf{U}_{|\mathcal{O}|})} \mathrm{C}_{|\mathcal{O}|}(k).
\end{equation}
Now let $\mathrm{Orb}(\mathcal{G}, \theta)$ be the set of $(\theta)$-orbits of $\mathcal{G}$. Since
$$
A = \bigoplus_{\mathcal{O} \in \mathrm{Orb}(\mathcal{G}, \theta)}  (\mathrm{sp}(\mathcal{O}) \times H)
$$
we have
\begin{equation}\label{EqSS2}
A = \bigoplus_{\mathcal{O} \in \mathrm{Orb}(\mathcal{G}, \theta)}\left(\bigoplus_{g \in R(G, \mathbf{U}_{|\mathcal{O}|})}  \mathrm{C}(\mathcal{O}, \mathbf{U}_{|\mathcal{O}|})(1 \times g)\right)
\end{equation}
as a coalgebra by (\ref{EqSS1}), and  $\mathrm{C}(\mathcal{O}, \mathbf{U}_{|\mathcal{O}|})(1 \times g) \simeq \mathrm{C}_{|\mathcal{O}|}(k)$. We have shown that $A$ is cosemisimple and have described its coalgebra structure explicitly. In particular
\begin{equation}\label{EqGAF}
\mathrm{G}(A) = \mathrm{F}(\theta, \mathcal{G}) \times \mathrm{G} \;\; \mbox{and thus} \;\;  |\mathrm{G}(A)| = |\mathrm{F}(\theta, \mathcal{G})||\mathrm{G}|,
\end{equation}
where $\mathrm{F}(\theta, \mathcal{G})$ is the subgroup of fixed points of $\theta$ in $\mathcal{G}$. These equations should be interesting in connection with the examples at the end of this section.

The Hopf algebra $A$ is involutory. Let $S$, $S_H$, and $S_B$ denote the antipodes of $A$, $H$, and $B$ respectively. Let $c \in B$ and $h \in H$. Then
\begin{equation}\label{EqAntipode}
S(c \times h) = (1 \times S_H(c_{(-1)}h))(S_B(c_{(0)}) \times 1)
\end{equation}
by part (b) of \cite[\S 2 Proposition 2]{RadProj}. Since $c \times h = (c \times 1)(1 \times h)$, to show that $S^2 = \mathrm{I}_A$ we need only show that the equation holds on $c \times 1$ and $1 \times h$. Since $S^2_H = \mathrm{I}_H$ the equation holds on the latter by (\ref{EqAntipode}). To show it holds on the former we may assume $0 \neq c \in B_g$ for some $g \in G$. Since $B_g \neq (0)$ it follows that $g{\cdot}c = c$. Since $S_B$ and $\theta$ commute, $S_B(B_{g}) \subseteq B_{g}$. We use (\ref{EqAntipode}) and the fact that $S_B^2 = \mathrm{I}_B$ to compute
\begin{eqnarray*}
S^2(c \times 1) & = & S((1 \times S_H(g))(S_B(c) \times 1)) \\
& = & S(S_B(c) \times 1)(1 \times g) \\
& = & (1 \times S_H(g))(S_B^2(c) \times 1)(1 \times g) \\
& = & (1 \times g^{-1})(c \times g) \\
& = & g^{-1}{\cdot}c \times g^{-1}g \\
& = & c \times 1.
\end{eqnarray*}
Therefore $S^2 = \mathrm{I}_A$.

We end this section with some examples. First we formally summarize the preceding discussion. To do this we introduce a bit of useful linear algebra and notation.

Let $T$ be an endomorphism of a vector space $V$. Then $T$ is \textit{diagonalizable} if $V$ is the sum of eigenspaces of $T$. Suppose $T$ is an automorphism and $V$ has a basis $B$ such that $T(B) = B$ and the $(T)$-orbits of $B$ are finite. Then $T$ is diagonalizable if and only if $k^\times$ contains a primitive $r^{th}$ root of unity for every orbit length $r$. If $\mathcal{G}$ is a group and $\theta \in \mathrm{Aut}_{\mathrm{Group}}(\mathcal{G})$ we let $\theta_{k[\mathcal{G}]} : k[\mathcal{G}] \mapsto k[\mathcal{G}]$ denote the linear extension of $\theta$ to $k[\mathcal{G}]$.

For a cardinal number $m$ and a vector space $V$ we denote the direct sum of $m$ copies of $V$ by $m{\cdot}V$.
\begin{Theorem}\label{ThmCoalgebraStructure}
Let $G$ and $\mathcal{G}$ be groups, let $\pi : G \longrightarrow \mathrm{Aut}_{\mathrm{Group}}(\mathcal{G})$ be a group homomorphism, and let $\theta \in \mathrm{Aut}_{\mathrm{Group}}(\mathcal{G})$. Suppose that the $(\theta)$-orbits of $\mathcal{G}$ are finite and that $\theta_{k[\mathcal{G}]}$ is diagonalizable. Assume further that:
\begin{enumerate}
\item[{\rm (a)}] $\pi(g){\circ}\theta = \theta{\circ}\pi(g)$ for all $g \in G$.
\item[{\rm (b)}] There exists an isomorphism $U \longrightarrow \mathbf{U}$ $(\lambda \mapsto \mbox{\boldmath $\lambda$})$ of the subgroup $U \subseteq k^\times$ generated by the eigenvalues of $\theta_{k[\mathcal{G}]}$ and a subgroup $\mathbf{U}$ of $G$ which satisfies $\mathbf{U} \subseteq \mathrm{Ker}(\pi)\cap \mathrm{Z}(G)$.
\end{enumerate}
Let $H = k[G]$ and $B = k[\mathcal{G}]$ Then:
\begin{enumerate}
\item[{\rm (c)}]  The Hopf algebra $B$ over $k$ is also a Hopf algebra of ${}^H_H\mathcal{YD}$, where its $H$-module structure is given by $g{\cdot}b = \pi(g)(b)$ for all $g \in G$ and $b \in \mathcal{G}$ and its $H$-comodule structure is given by $\rho(b) = \mbox{\boldmath $\lambda$} \otimes b$ for $b \in B_{(\lambda)}$.
\item[{\rm (d)}] The biproduct $B \times H$ is a cosemisimple involutory Hopf algebra over $k$. In particular
$$
B \times H \simeq \bigoplus_{\mathcal{O} \in \mathrm{Orb}(\mathcal{G}, \theta)}[G :  \mathbf{U}_{|\mathcal{O}|}]{\cdot}\mathrm{C}_{|\mathcal{O}|}(k)
$$
as coalgebras.
\end{enumerate}
\qed
\end{Theorem}

We present a litany of examples. Non-trivial semisimple cosemisimple Hopf algebras of dimension a power of $2$ over certain fields $k$ have been the subject of several papers by Kashina \cite{KashYevg2000Dim16,KashYevg2003SS2m,KashYevg2003SS2fam2M,KashYevg2007}. The reader should compare these with the ones constructed below whose dimension is a power of $2$.

\begin{Ex}\label{Ex1}
$\mathrm{A}(\mathcal{G}, \theta) = B \times H$, where $B = k[\mathcal{G}]$ and $H =  k[G]$, constructed as above: $\mathcal{G}$ is any group; $\theta \in \mathrm{Aut}_{\mathrm{Group}}(\mathcal{G})$, $\theta_{k[\mathcal{G}]}$ is diagonalizable, and the $(\theta)$-orbits of $\mathcal{G}$ are finite; $G = U \times (\theta)$, where $U$ is the subgroup of $k^\times$ generated by the eigenvalues of $\theta$;  $\pi : G \longrightarrow \mathrm{Aut}_{\mathrm{Group}}(\mathcal{G})$ is given by the projection of $G$ onto $(\theta)$; $\mathbf{U} = \mathrm{Ker}(\pi)$, the first factor of $G$; and $U \simeq \mathbf{U}$ is the usual identification.
\end{Ex}

Comments on Example \ref{Ex1}. The identification of (\ref{EqGAF}) is also one of groups. Note (\ref{EqCommute}) holds, since $\mathrm{Im}(\pi) = (\theta)$, and $\mathbf{U} \subseteq \mathrm{Z}(G)$ since $G$ is commutative. Suppose $\theta \neq \mathrm{I}_\mathcal{G}$. Then $\mathrm{A}(\mathcal{G}, \theta)$ not a group algebra by (\ref{EqSS2}) and is not commutative since $\mathrm{Ker}(\pi) \neq G$. In this case when $\mathcal{G}$ and $G$ are finite, $A$ is a finite-dimensional non-trivial semisimple cosemisimple Hopf algebra.

\begin{Ex}\label{Ex2}
A special case of Example \ref{Ex1}, where the characteristic of $k$ is not $2$, $\mathcal{G} = \mathbb{Z}_n$, $n \geq 3$ and odd, $\theta(x) = -x$ for all $x \in \mathcal{G}$. Here $H \simeq k[\mathbb{Z}_2 \times \mathbb{Z}_2]$. Then $A = \mathrm{A}(\mathcal{G}, \theta)$ is a non-trivial semisimple cosemisimple Hopf algebra over $k$, $\mathrm{Dim}(A) = 4n$, $\mathrm{G}(A) \simeq \mathbb{Z}_2 \times \mathbb{Z}_2$, and $\mathrm{Rank}_{k[\mathrm{G}(A)]}A = n$.
\end{Ex}

\begin{Ex}\label{Ex2A}
Example \ref{Ex2} with the difference $n = 2m \geq 3$. Again $A = \mathrm{A}(\mathcal{G}, \theta)$ is a non-trivial semisimple cosemisimple Hopf algebra over $k$, $\mathrm{Dim}(A) = 4n$, $\mathrm{G}(A) \simeq \mathbb{Z}_2 \times \mathbb{Z}_2 \times \mathbb{Z}_2$, and $\mathrm{Rank}_{k[\mathrm{G}(A)]}A = m$.
Observe that the preceding example accounts for a non-trivial semisimple cosemisimple Hopf algebra $A$ of dimension $2^n$, where $n \geq 4$ and the characteristic of $k$ is not $2$, such that $\mathrm{Rank}_{k[\mathrm{G}(A)]}A = 2^{n - 3}$.
\end{Ex}

\begin{Ex}\label{Ex2B}
A special case of Example \ref{Ex1}, where the characteristic of $k$ is not $2$, $\mathcal{G} = \mathbb{Z}_n \times \mathbb{Z}_n$, $n \geq 2$, $\theta((x, y)) = (y, x)$ for all $(x, y) \in \mathcal{G}$. Here $H \simeq k[\mathbb{Z}_2 \times \mathbb{Z}_2]$. Then $A = \mathrm{A}(\mathcal{G}, \theta)$ is a non-trivial semisimple cosemisimple Hopf algebra over $k$, $\mathrm{Dim}(A) = 4n^2$, $\mathrm{G}(A) \simeq \mathbb{Z}_n \times \mathbb{Z}_2 \times \mathbb{Z}_2$, and $\mathrm{Rank}_{k[\mathrm{G}(A)]}A = n$.
\end{Ex}

\begin{Ex}\label{Ex6}
A special case of Example \ref{Ex1} with $\mathcal{G} = \mathbb{Z}_2 \times \mathbb{Z}_2$. Let $\theta \in \mathrm{Aut}_\mathrm{Group}(\mathcal{G}) = S_3$ be a $3$-cycle. Assume that $k$ contains a primitive $3^{rd}$ root of unity. Then $\theta_{k[\mathcal{G}]}$ is diagonalizable. $H \simeq k[\mathbb{Z}_3 \times \mathbb{Z}_3]$ and $A = \mathrm{A}(\mathcal{G}, \theta)$ is a non-trivial semisimple cosemisimple Hopf algebra of dimension $36$. Observe that $\mathrm{G}(A) \simeq \mathbb{Z}_3 \times \mathbb{Z}_3$.
\end{Ex}

Let $m, n \geq 2$. Given a $2$-group $G$ of order $2^n$ there exists a non-trivial semisimple cosemisimple Hopf algebra $A$ over $k$ of dimension $2^{m + n}$ with a Hopf algebra projection onto $H = k[G]$. Thus $\mathrm{rank}_H(A) = 2^m$. Compare with Corollary \ref{Cor2}.
\begin{Ex}\label{Ex4}
Let $\mathcal{G}$ and $G$ be $2$-groups of orders $2^m$ and $2^n$ respectively, where $m, n \geq 2$, and suppose that $\mathcal{G}$ is abelian. Suppose that the characteristic of $k$ is not $2$. Since $\mathcal{G}$ is abelian of order at least $4$, there exists $\theta \in \mathrm{Aut}_\mathrm{Group}(\mathcal{G})$ of order $2$. Since $G$ is a $2$-group of order $2^n$, and $n \geq 2$ it follows that $G$ contain a normal subgroup $N$ of index $2$ which contains a $2$-element subgroup $\mathbf{U} \subseteq \mathrm{Z}(G)$. Since $(\theta)$ has order $2$ there exists a group homomorphism $\pi : G \longrightarrow \mathrm{Aut}_\mathrm{Group}(\mathcal{G})$ such that $\mathrm{Im}(\pi) = (\theta)$ and $\mathrm{Ker}(\pi) = N$. Identify $\mathbf{U}$ with the group $\{-1, 1 \}$ of eigenvalues of $\theta$.

Let $B = k[\mathcal{G}]$ and $H = k[G]$. The resulting biproduct $A = B \times H$ is a non-trivial semisimple cosemisimple Hopf algebra of dimension $2^{m + n}$. In particular $\mathrm{rank}_H(A) = 2^m$.
\end{Ex}

For the assertions concerning subgroups in the previous example see the proof of \cite[Proposition 2.106]{Rotman2002} for example.

We recall the example described on page 342 of \cite{RadProj} which was constructed by the methods of this section.
\begin{Ex}\label{Ex3}
Let $\mathcal{G}$ be any commutative group and $\theta \in \mathrm{Aut}_\mathrm{Group}(\mathcal{G})$ have order $n \geq 1$, $G = (a)$ cyclic of order $2n$. Assume $k$ has a primitive $n^{th}$ root of unity. $\pi : G \longrightarrow \mathrm{Aut}_\mathrm{Group}(\mathcal{G})$ is determined by $\pi (a) = S$, where $S(b) = b^{-1}$ for all $b \in \mathcal{G}$. Identify the subgroup $\mathrm{Ker}(\pi) = (a^2)$, which has order $n$, with the group of $n^{th}$ roots of unity of $k^\times$, the set of eigenvalues of $\theta$. The resulting biproduct $A = k[\mathcal{G}] \times k[G]$ is cosemisimple and not a group algebra, $n \geq 2$.
\end{Ex}

Apropos of the preceding example, suppose that $\mathcal{G}$ is finite of order $m$ and $S, \theta \neq \mathrm{I}_{\mathcal{G}}$. Then $A$ is a non-trivial semisimple cosemisimple Hopf algebra over $k$, $\mathrm{Dim}(A) = 2mn$. Such an example was constructed on page 342 of \cite{RadProj} as one of a family. Here $m = 3$ and $n = 2$. Thus $\mathrm{Dim}(A) = 12$ and $\mathrm{G}(A) = \mathbb{Z}_4$. There are only two non-trivial $12$-dimensional semisimple cosemisimple Hopf algebras over an algebraically closed field $k$ of characteristic not $2$ or $3$. This was established in \cite{Fukuda}. We note here that the second one is accounted for by Example \ref{Ex2} with $\mathcal{G} = \mathbb{Z}_3$ also. Here $\mathrm{G}(A) = \mathbb{Z}_2 \times \mathbb{Z}_2$; thus the two $12$-dimensional Hopf algebras described here can not be isomorphic.

Another $36$-dimensional Hopf algebra is obtained by virtue of Example \ref{Ex3}. Assume that $k$ contains a primitive $3^{rd}$ root of unity. We start with the group $\mathcal{G} = \mathbb{Z}_2 \times \mathbb{Z}_2$ and automorphism $\theta$ of Example \ref{Ex6} and let $G = (a)$ to be cyclic of order $9$ and $\pi$ be determined by $\pi(a) = \theta$. In this case $\mathrm{G}(A) \simeq \mathbb{Z}_9$.

The algebra structure of our examples has the form $k[\mathcal{G}] \# k[G]$ for certain groups $\mathcal{G}, G$. With the exception of Example \ref{Ex4} the group $G$ is a finite abelian. We are interested in finite-dimensional examples. Here, with the exception of the generic Example \ref{Ex1}, the group $\mathcal{G}$ is abelian, and hence $k[\mathcal{G}]$ is a semisimple algebra when $k$ is algebraically closed of characteristic zero. In this case $k[\mathcal{G}]$ is the product of copies of $k$ as a $k$-algebra. We begin our discussion of the algebra $k[\mathcal{G}] \# k[G]$ from a slightly more general perspective, and will assume $k$ has ample roots of unity.

\section{The smash product $B \# k[G]$, where $G$ is a finite abelian group}

First some comments on smash products. Let $H$ and $L$ be bialgebras over the field $k$ and supposed that $B$ is a left $H \otimes L$-module algebra. Then $B$ is a left $H$-module and a left $L$-module via pullback along $H \longrightarrow H \otimes L$ $(h \mapsto h \otimes 1)$ and $L \longrightarrow H \otimes L$ $(\ell \mapsto 1 \otimes \ell)$ respectively. With the pullback actions, which commute, $B$ is a left $H$-module and a left $L$-module algebra. Observe that $B \# H$ is a left $L$-module algebra, where $\ell{\cdot}(b \# h) = \ell{\cdot}b \# h$, and $(B \# H) \# L \simeq B \# (H \otimes L)$ as $k$-algebras, where $(b \# h) \# \ell \mapsto b \# (h \otimes \ell)$. Conversely, suppose that $B$ is a left $H$-module algebra and a left $L$-module algebra and that the actions commute. Then $B$ is a left $H \otimes L$-module algebra, where $(h \otimes \ell){\cdot}b = h{\cdot}(\ell{\cdot}b)$. The actions on $H$ and $L$ are obtained by pullback as above and we have the above isomorphism of smash products.

Special cases of Example \ref{Ex1} account for most of the Hopf algebras in this paper. By virtue of the discussion above we have
$$
\mathrm{A}(\mathcal{G}, \theta) = k[\mathcal{G}] \# k[U \times (\theta)] \simeq k[\mathcal{G}] \# (k[U] \otimes k[(\theta)]) \simeq (k[\mathcal{G}] \# k[U]) \# k[(\theta)].
$$
When $\mathcal{G}$ is abelian $k[\mathcal{G}] \# k[U] \simeq k[\mathcal{G}] \otimes k[U]$ and therefore
\begin{equation}\label{EqAGTheta}
\mathrm{A}(\mathcal{G}, \theta) \simeq (k[\mathcal{G}] \otimes k[U]) \# k[(\theta)]
\end{equation}
as an algebra over $k$.

To study the ideal structure of the $\mathrm{A}(\mathcal{G}, \theta)$'s we need only consider smash products of the form $B \# k[G]$, where $G$ is finite cyclic. We will take $G$ to be a finite abelian group. Since $k[G]$ is the tensor product of group algebras of cyclic groups, we could reduce to the cyclic case. We will find it a bit more interesting not to do this.

Let $n \geq 2$ and suppose the field and $k$ has a primitive $n^{th}$ root of unity $\lambda$. To set notation, first let $G = (g)$ be a cyclic group of order $n$. For $i \geq 0$ we define
\begin{equation}\label{EqeiOrthogonal}
e_i = \sum_{\ell = 0}^{n-1} (\lambda^{i\ell}/n)g^\ell.
\end{equation}
Observe that if $i \equiv i' \; (\mathrm{mod} \, n)$ then $\lambda^i = \lambda^{i'}$, $g^i = g^{i'}$, and therefore $e_i = e_{i'}$. We can regard scripts as elements of $\mathbb{Z}_n$. Thus $e_0, \ldots, e_{n-1}$ lists the distinct $e_i$'s and these elements can be described by
\begin{equation}\label{EqEI}
e_i = \sum_{\ell \in \mathbb{Z}_n} (\lambda^{i\ell}/n)g^\ell
\end{equation}
for $i \in \mathbb{Z}_n$. Notice that $\{e_0, \ldots, e_{n-1}\} = \{e_i\}_{i \in \mathbb{Z}_n}$ is an orthogonal set of idempotents which is a linear basis for $H = k[G]$. Thus
$$
\sum_{i \in \mathbb{Z}_n}e_i = 1 \;\; \mbox{and} \;\; e_ie_j = \delta_{i, j}e_i, \;\; \mbox{and we note} \;\; g^j e_i = \lambda^{-ij}e_i,
$$
for all $i, j \in \mathbb{Z}_n$. In particular $H$ is a semisimple algebra whose minimal ideals are the one-dimensional subspaces $ke_i$.

The coalgebra structure of $H$ in terms of this basis is given by
$$
\Delta(e_i) = \sum_{\ell \in \mathbb{Z}_n} e_\ell \otimes e_{i - \ell} \;\; \mbox{and} \;\; \epsilon(e_{i}) = \delta_{i, 0}
$$
for all $i \in \mathbb{Z}_n$. Note that the unique integral $\Lambda$ for $H$ which satisfies $\epsilon(\Lambda) = 1$ is $\Lambda = e_0$.

Now let $G = G_1 \times \cdots \times G_t$ be a product of cyclic groups, where $G_j = (g_j)$ and has order $n_j$ for all $1 \leq j \leq t$. Suppose $G$ has order $n$. For each $1 \leq j \leq t$ observe that $\lambda_j = \lambda^{n/n_j}$ is a primitive $n_j^{th}$ root of unity.

Let $\mathbf{G} = \mathbb{Z}_{n_1} \oplus \cdots \oplus \mathbb{Z}_{n_t}$ with its usual \textit{ring} structure. We define formal powers. For $\mathbf{m} = (m_1, \ldots, m_t) \in \mathbf{G}$ set
$$
g^{(\mathbf{m})} = (g_1^{m_1}, \ldots, g_t^{m_t}) \;\; \mbox{and} \;\; \lambda^{(\mathbf{m})} = \lambda_1^{m_1} \cdots \lambda_t^{m_t}.
$$
Note that $g^{(-\mathbf{m})} = (g^{(\mathbf{m})})^{-1}$ and $\lambda^{(-\mathbf{m})} = (\lambda^{(\mathbf{m})})^{-1}$ for $\mathbf{m} \in \mathbf{G}$. Observe that $g^{(\mathbf{m})}g^{(\mathbf{n})} = g^{(\mathbf{m} + \mathbf{n})}$ and $\lambda^{(\mathbf{m})}\lambda^{(\mathbf{n})} = \lambda^{(\mathbf{m} + \mathbf{n})}$ for all $\mathbf{m}, \mathbf{n} \in \mathbf{G}$. The map $\mathbf{G} \longrightarrow G$ defined by $\mathbf{m} \mapsto g^{(\mathbf{m})}$ is a group isomorphism.

A few more definitions. For $1 \leq j \leq t$ define $e(j)_i$ as in (\ref{EqEI}), where $\mathbb{Z}_{n_j}$ replaces $\mathbb{Z}_n$, $\lambda_j$ replaces $\lambda$, and $g_j$ replaces $g$. For $\mathbf{m} = (m_1, \ldots, m_t) \in \mathbf{G}$ let $e_{\mathbf{m}}$ be the image of $e(1)_{m_1} \otimes \cdots \otimes e(t)_{m_t}$ under the algebra isomorphism $k[G_1] \otimes \cdots \otimes k[G_t] \simeq k[G_1 \times \cdots \times G_t] = k[G]$ determined by $g_1^{m_1} \otimes \cdots \otimes g_t^{m_t} \mapsto (g_1^{m_1}, \ldots, g_t^{m_t}) = g^{(\mathbf{m})}$.
Our formulas for the $e_i$'s defined in the cyclic case have generalizations to $G = G_1 \times \cdots \times G_t$ which the reader can easily verify. First of all note that
$$
e_{\mathbf{m}}= \sum_{\mathbf{r} \in {\mathbf{G}}} (\lambda^{(\mathbf{m}\mathbf{r})}/n)g^{(\mathbf{r})}
$$
for all $\mathbf{m} \in \mathrm{G}$ and $\{e_{\mathbf{m}}\}_{\mathbf{m} \in \mathbf{G}}$ is an orthogonal set of idempotents which is a linear basis for $H = k[G]$. Thus
$$
\sum_{\mathbf{m} \in \mathbf{G}} e_{\mathbf{m}} = 1 \;\; \mbox{and} \;\; e_{\mathbf{m}} e_{\mathbf{n}} = \delta_{\mathbf{m}, \mathbf{n}}, \;\; \mbox{and we also note} \;\; g^{(\mathbf{m})}e_{\mathbf{n}} = \lambda^{(-\mathbf{m}\mathbf{n})}e_{\mathbf{n}},
$$
for all $\mathbf{m}, \mathbf{n} \in \mathbf{G}$. As a consequence $H$ is a semisimple algebra whose minimal ideals are the one-dimensional subspaces $ke_{\mathbf{m}}$.

The coalgebra structure of $H$ in terms of this basis is given by
$$
\Delta(e_\mathbf{m}) = \sum_{\mathbf{j} \in \mathbf{G}} e_\mathbf{j} \otimes e_{\mathbf{m} - \mathbf{j}} \;\; \mbox{and} \;\; \epsilon(e_\mathbf{m}) = \delta_{\mathbf{m}, \mathbf{0}}
$$
for all $\mathbf{m} \in \mathbf{G}$, where $\mathbf{0}$ is the neutral element for the additive group structure of $\mathbf{G}$. Note that the unique integral $\Lambda$ for $H$ which satisfies $\epsilon(\Lambda) = 1$ is $\Lambda = e_\mathbf{0}$.

Now let $B$ be a left $H$-module algebra. Then for $\mathbf{m} \in \mathbf{G}$ the linear automorphism $\theta_\mathbf{m}$ of $B$ defined by $\theta_\mathbf{m}(x) = g^{(\mathbf{m})}{\cdot}x$ for all $x \in B$ is an algebra automorphism since $g^{(\mathbf{m})} \in \mathrm{G}(H)$. We will use the fact that
\begin{equation}\label{Eqeibejc}
(e_\mathbf{m}{\cdot}b)(e_\mathbf{n}{\cdot}c) = e_{\mathbf{m} + \mathbf{n}}{\cdot}(b(e_\mathbf{n}{\cdot}c))
\end{equation}
for all $\mathbf{m}, \mathbf{n} \in \mathbf{G}$ and $b, c \in B$. To see this we calculate
\begin{eqnarray*}
(e_\mathbf{m}{\cdot}b)(e_\mathbf{n}{\cdot}c) & = &
\sum_{\mathbf{p}, \mathbf{q} \in {\mathbf{G}}} \left(((\lambda^{(\mathbf{m}\mathbf{p})}/n)g^{(\mathbf{p})}{\cdot}b)((\lambda^{(\mathbf{n}\mathbf{q})}/n)g^{(\mathbf{q})}{\cdot}c)\right) \\
& = & \sum_{\mathbf{p}, \mathbf{q} \in {\mathbf{G}}}(\lambda^{(\mathbf{m}\mathbf{p} + \mathbf{n}\mathbf{q})}/n^2)g^{(\mathbf{p})}{\cdot}(b(g^{(\mathbf{q} - \mathbf{p})}{\cdot}c)) \\
& = & \sum_{\mathbf{p}, \mathbf{r} \in {\mathbf{G}}}(\lambda^{(\mathbf{m}\mathbf{p} + \mathbf{n}(\mathbf{r} + \mathbf{p}))}/n^2)g^{(\mathbf{p})}{\cdot}(b(g^{(\mathbf{r})}{\cdot}c)) \\
& = & \sum_{\mathbf{p}, \mathbf{r} \in {\mathbf{G}}}(\lambda^{((\mathbf{m} + \mathbf{n})\mathbf{p}))}/n)(\lambda^{(\mathbf{n}\mathbf{r})}/n)g^{(\mathbf{p})}{\cdot}(b(g^{(\mathbf{r})}{\cdot}c)) \\
& = & e_{\mathbf{m} + \mathbf{n}}{\cdot}(b(e_\mathbf{n}{\cdot}c)).
\end{eqnarray*}

Let $A = B \# H$ be the smash product of $B$ and $H$. Recall multiplication is defined by $(b \# h)(b' \# h') = b(h_{(1)}{\cdot}b') \# h_{(2)}h'$ for all $b, b' \in B$ and $h, h' \in H$ and $1_A = 1_B \# 1_H$. Again, $(b \# h)(b' \# h') = bb' \# hh'$ if $b' = 1$ or $h = 1$.

The ideals of $A$ are easy to characterize. First note that
$$
1_A = 1_B \# 1_H = \sum_{\mathbf{m} \in \mathbf{G}} 1 \# e_{\mathbf{m}}
$$
and $(1 \# e_{\mathbf{m}})(1 \# e_{\mathbf{n}}) = \delta_{{\mathbf{m}}, {\mathbf{n}}}(1 \# e_{\mathbf{m}})$ for all ${\mathbf{m}}, {\mathbf{n}} \in {\mathbf{G}}$. Let $I$ be an ideal of $A$. Then $I(1 \# e_\mathbf{m}) \subseteq I$ for all $\mathbf{m} \in \mathbf{G}$ and thus
$$
I = \bigoplus_{\mathbf{m} \in \mathbf{G}} I(1 \# e_\mathbf{m})
$$
and
since $I$ is a right ideal of $A$. For $\mathbf{m}$ we define $I_\mathbf{m} = \{b \in B \, | \, b \# e_\mathbf{m} \in I\}$. Since $(b \# 1)(b' \# e_\mathbf{m}) = bb' \# e_\mathbf{m}$ and $(1 \# g^{(\mathbf{r})})(b' \# e_\mathbf{m}) = g^{(\mathbf{r})}{\cdot}b' \# \lambda^{(-\mathbf{r}\mathbf{m})}e_\mathbf{m}$ for all $b, b' \in B$, it follows that $I_\mathbf{m}$ is a left ideal and a left $H$-submodule of $B$ since $I$ is a left ideal of $A$.

We show that $I_\mathbf{m} \# e_\mathbf{m} = I(1 \# e_\mathbf{m})$. By definition
$$
I_\mathbf{m} \# e_\mathbf{m} = (I_\mathbf{m} \# e_\mathbf{m})(1 \# e_\mathbf{m}) \subseteq I(1 \# e_\mathbf{m}).
$$
On the other hand let $a \in I(1 \# e_\mathbf{m})$. Since $\{e_\mathbf{r}\}_{\mathbf{r} \in \mathbf{G}}$ is a basis for $H$ we may write $a = \sum_{\mathbf{r} \in \mathbf{G}} b_\mathbf{r} \# e_\mathbf{r}$ where $b_\mathbf{r} \in B$. Since $a = a(1 \# e_\mathbf{m}) = b_\mathbf{m} \# e_\mathbf{m}$, by definition $b_\mathbf{m} \in I_\mathbf{m}$. Therefore $a \in I_\mathbf{m} \# e_\mathbf{m}$. We have shown $I_\mathbf{m} \# e_\mathbf{m} = I(1 \# e_\mathbf{m})$.

Let $b, b' \in B$, and $\mathbf{m}, \mathbf{r} \in \mathbf{G}$. We observe that
\begin{equation}\label{EqMain}
b(e_{\mathbf{m} - \mathbf{r}}{\cdot}b') \# e_\mathbf{r} = (b \# e_\mathbf{m})(b' \# e_\mathbf{r}) = (b \# e_\mathbf{m})(e_{\mathbf{m} - \mathbf{r}}{\cdot} b' \# 1).
\end{equation}

Now suppose $b \in I_\mathbf{m}$ in the first equation of (\ref{EqMain}). Then $b(e_{\mathbf{m} - \mathbf{r}}{\cdot}b') \in I_\mathbf{r}$ since the right hand side of the equation belongs to $I$. We have established part (a) of the following:
\begin{Prop}\label{PropIdealsOfA}
Let $G$ be the abelian group of order $n$ described above, let $k$ be any field which contains a primitive $n^{th}$ root of unity, let $H = k[G]$, and let the $e_\mathbf{m}$'s be as above. Let $B$ be a left $H$-module algebra and $A = B \# H$.
\begin{enumerate}
\item[{\rm (a)}] If $I$ is an ideal of $A$ then $I = \bigoplus_{\mathbf{m} \in \mathbf{G}} I_\mathbf{m} \# e_\mathbf{m}$, where (i) each $I_\mathbf{m}$ is a left ideal and a left $H$-submodule of $B$ and (ii) $I_\mathbf{m}(e_{\mathbf{m} - \mathbf{r}}{\cdot}B) \subseteq I_\mathbf{r}$ for all $\mathbf{r} \in \mathbf{G}$.
\item[{\rm (b)}] Suppose $\{ I_\mathbf{m}\}_{\mathbf{m} \in \mathbf{G}}$ are subspaces of $B$ which satisfy the conditions of part (a). Then $I = \bigoplus_{\mathbf{m} \in \mathbf{G}} I_\mathbf{m} \# e_\mathbf{m}$ is an ideal of $A$.
\end{enumerate}
\end{Prop}

\pf
We need only establish part (b). Suppose the hypothesis of part (b) holds and let $I = \bigoplus_{\mathbf{m} \in \mathbf{G}} I \# e_\mathbf{m}$. For $\mathbf{m} \in \mathbf{G}$ we see by the first equation of (\ref{EqMain}) and (i) that
$$
(B \# e_\mathbf{r})(I_\mathbf{m} \# e_\mathbf{m}) \subseteq B(e_{\mathbf{r} - \mathbf{m}} {\cdot}I_\mathbf{m}) \# e_\mathbf{m} \subseteq I_\mathbf{m} \# e_\mathbf{m} \subseteq I
$$
and by the same and (ii) that
$$
(I_\mathbf{m} \# e_\mathbf{m})(B \# e_\mathbf{r}) \subseteq  I_\mathbf{m}(e_{\mathbf{m} - \mathbf{r}}{\cdot}B) \# e_\mathbf{r} \subseteq  I_\mathbf{r} \# e_\mathbf{r} \subseteq I.
$$
Since the $B \# e_\mathbf{r}$'s span $A$, it follows that $I$ is an ideal of $A$.
\qed
\medskip

The first equation of (\ref{EqMain}) establishes the initial assertion of:
\begin{Cor}\label{CorLeftIdealA}
Assume the hypothesis of the preceding proposition and that $L$ is a left ideal and a left $H$-submodule of $B$. Then $L \# e_\mathbf{m}$ is a left ideal of $A$ for all $\mathbf{m} \in \mathbf{G}$. Furthermore the ideal of $A$ it generates is $(L \# e_\mathbf{m})(B \# 1)$.
\end{Cor}

\pf
Fix $\mathbf{m} \in \mathbf{G}$. Since $L \# e_\mathbf{m}$ is a left ideal of $A$, $I = (L \# e_\mathbf{m})(B \# 1)$ is a left ideal of $A$ also. For $\mathbf{r} \in \mathbf{G}$ note that
$$
I(B \# e_\mathbf{r}) = (L \# e_\mathbf{m})(B \# 1)(B \# e_\mathbf{r}) = (L \# e_\mathbf{m})(B \# e_\mathbf{r}).
$$
Thus by the second equation of (\ref{EqMain}) we have
$$
I(B \# e_\mathbf{r}) = (L \# e_\mathbf{m})(e_{\mathbf{m} - \mathbf{r}}{\cdot}B \# 1) \subseteq I.
$$
Again, since the $B \# e_\mathbf{r}$'s span $A$, it follows that $I$ is a right ideal of $A$ also.
\qed
\medskip

We will say that a left $H$-module algebra $B$ is \textit{semisimple as a module algebra} if whenever $L$ is a left ideal and a left $H$-submodule of $B$ then $L \oplus L' = B$ for some left ideal and left $H$-submodule $L'$ of $B$.
\begin{Cor}\label{CorLeftIdealModA}
Assume the hypothesis of the preceding proposition and in addition that $B$ is semisimple as a module algebra. Let $L$ be an ideal and a left $H$-submodule of $B$. Then the ideal of $A$ generated by $L \# e_\mathbf{m}$ is $(L \# e_\mathbf{m})(L \# 1)$.
\end{Cor}

\pf
By assumption $L \oplus L' = B$ for some left ideal and left $H$-submodule $L'$ of $B$. Since $L$ is an ideal of $B$, $LL + LL' = L$; thus $LL'= (0)$ and $LL = L$. Now the ideal of $A$ which $(L \# e_\mathbf{m})$ generates is $I = (L \# e_\mathbf{m})(B \# 1)$ by the preceding corollary. Let $\mathbf{n} \in \mathbf{G}$. Then
$$
(L \# e_\mathbf{m})(L' \# e_\mathbf{n}) \subseteq L(e_{\mathbf{m}- \mathbf{n}}{\cdot} L') \# e_\mathbf{n} \subseteq LL' \# e_\mathbf{n} = (0)
$$
by the first equation of (\ref{EqMain}). Since $L' \# 1 \subseteq \sum_{\mathbf{n} \in \mathbf{G}} L' \# e_\mathbf{n}$ it follows that $(L \# e_\mathbf{m})(L' \# 1) = (0)$. As a consequence
$$
(L \# e_\mathbf{m})(B \# 1) = (L \# e_\mathbf{m})(L \# 1) + (L \# e_\mathbf{m})(L' \# 1) = (L \# e_\mathbf{m})(L \# 1).
$$
Corollary \ref{CorLeftIdealA} completes the proof.
\qed
\medskip

We end this section with a discussion of the minimal ideals of the smash product $A$ of Proposition \ref{PropIdealsOfA}. Suppose first of all that $I$ is any non-zero ideal of $A$. Since $I = \bigoplus_{\mathbf{m} \in \mathbf{G}} I_\mathbf{m} \# e_\mathbf{m}$ it follows that $I_{\mathbf{m}_0} \neq (0)$ for some $\mathbf{m}_0 \in \mathbf{G}$. Let $S$ be a non-zero left ideal and $H$-submodule of $B$ contained in $I_{\mathbf{m}_0}$. Define $S_\mathbf{m} = S(e_{\mathbf{m}_0 - \mathbf{m}}{\cdot}B)$ for $\mathbf{m} \in \mathbf{G}$. Then each $S_\mathbf{m}$ is a left ideal and a left $H$-submodule of $B$. Observe that some $S_\mathbf{m} \neq (0)$ for some $\mathbf{m}$ since $S \subseteq SB = \sum_{\mathbf{m} \in \mathbf{G}}S(e_\mathbf{m}{\cdot}B)$.

Let $\mathbf{m} \in \mathbf{G}$. Then $S_\mathbf{m} \subseteq L_\mathbf{m}$ since
$$
S_\mathbf{m} = S(e_{\mathbf{m}_0 - \mathbf{m}}{\cdot}B) \subseteq I_{\mathbf{m}_0}(e_{\mathbf{m}_0 - \mathbf{m}}{\cdot}B) \subseteq I_\mathbf{m}.
$$
For $\mathbf{m}, \mathbf{r} \in \mathbf{G}$ we have $S_\mathbf{m}(e_{\mathbf{m} - \mathbf{r}}{\cdot}B) \subseteq S_\mathbf{r}$ since
$$
S_\mathbf{m}(e_{\mathbf{m} - \mathbf{r}}{\cdot}B) = S(e_{\mathbf{m}_0 - \mathbf{m}}{\cdot}B)(e_{\mathbf{m} - \mathbf{r}}{\cdot}B) \subseteq S(e_{\mathbf{m}_0 - \mathbf{r}}{\cdot}(B(e_{\mathbf{m} - \mathbf{r}}{\cdot}B))) \subseteq S_\mathbf{r}
$$
by (\ref{Eqeibejc}). Therefore $J = \bigoplus_{\mathbf{m} \in \mathbf{G}}S_\mathbf{m} \# e_\mathbf{m}$ is a non-zero ideal of $A$ and $J \subseteq I$.

Suppose further that $I$ is minimal. Then $J = I$; in particular $I_\mathbf{m} = S_\mathbf{m} = S(e_{\mathbf{m}_0 - \mathbf{m}}{\cdot}B)$ for all $\mathbf{m} \in \mathbf{G}$. In particular
$$
(0) \neq I_{\mathbf{m}_\mathbf{0}} = S(e_\mathbf{0}{\cdot}B) = S(\Lambda{\cdot}B) = \sum_{b \in B}S(\Lambda{\cdot}b).
$$
Let $b \in B$. Then $S(\Lambda{\cdot}b)$ is a left ideal and a left $H$-submodule of $B$. Observe that $f_b : S \longrightarrow S(\Lambda{\cdot}b)$ defined by $f_b(s) = s(\Lambda{\cdot}b)$ for all $s \in S$ is a map of left $B$-modules and left $H$-modules. Thus if $S$ contains no proper subspace which is a left ideal and a left $H$-submodule of $B$ then $S(\Lambda{\cdot}b) = (0)$ or $f_b$ is an isomorphism.

\section{The case $B$ is a power of $k$}\label{SecLast}
We refer to the $k$-algebra $k \times  \cdots \times k$ with $m \geq 1$ factors as a \textit{power of} $k$. We continue the previous section and assume that $B$ a power of $k$. Then $B$ has has a basis $\mathcal{F}$ of idempotents which satisfy $1 = \sum_{f \in \mathcal{F}}f$ and $ff' = \delta_{f, f'}$ for all $f, f' \in \mathcal{F}$. Thus $B = \bigoplus_{f \in \mathcal{F}}kf$, which is also the direct sum of its minimal ideals, $B$ is commutative, and $B$ is semisimple. We note that $k[\mathcal{G}]$ of Example \ref{Ex1} is such an algebra $B$ when $\mathcal{G}$ is a finite abelian group and $k$ contains a primitive $|\mathcal{G}|^{th}$ root of unity.

Since the $\theta_\mathbf{m}$'s are algebra automorphisms of $B$ they permute the minimal ideals of $B$. Thus $\theta_\mathbf{m}(\mathcal{F}) = \mathcal{F}$ for all $\mathbf{m} \in \mathbf{G}$, which is to say that $G$ acts on $\mathcal{F}$ via the $H$-module action. For $f \in \mathcal{F}$ let $\mathcal{O}_f$ be the $G$-orbit of $f$ in $\mathcal{F}$. Observe that $L$ is a left ideal and a left $H$-submodule of $B$ if and only if $L = \mathrm{sp}(\mathcal{O}_{f_1}) \oplus \cdots \oplus \mathrm{sp}(\mathcal{O}_{f_s})$ is the direct sum of spans of orbits, where $f_1, \ldots, f_s \in \mathcal{F}$.

For $f \in \mathcal{F}$ and $\mathbf{m} \in \mathbf{G}$ let $\mathcal{M}_{f, \mathbf{m}} = (\mathrm{sp}(\mathcal{O}_f) \# e_\mathbf{m})(B \# 1)$. Then $\mathcal{M}_{f, \mathbf{m}}$ is an ideal of $A$ by Corollary \ref{CorLeftIdealA}. We will show that the $\mathcal{M}_{f, \mathbf{m}}$'s are the minimal ideals of $A$ and find necessary and sufficient conditions for two to be the same.

First an estimate of the dimension of $\mathcal{M}_{f, \mathbf{m}}$. Let $r = |\mathcal{O}_f|$ and let $L = \mathrm{sp}(\mathcal{O}_f)$. Then $L$ is an ideal and a left $H$-submodule of $B$. Since $B$ is semisimple as a module algebra
\begin{equation}\label{EqMfm}
\mathcal{M}_{f, \mathbf{m}} = (L \# e_\mathbf{m})(L \# 1)
\end{equation}
by Corollary \ref{CorLeftIdealModA}. Since $\mathrm{Dim}(L) = r$ it follows that $\mathrm{Dim}(\mathcal{M}_{f, \mathbf{m}}) \leq r^2$.

We next show that there is a linear isomorphism $\mathcal{M}_{f, \mathbf{m}} \simeq \mathrm{M}_r(k)$ of multiplicative structures. This will imply that $\mathcal{M}_{f, \mathbf{m}}$ is a minimal ideal of $A$. Choose an $r$-element subset $S \subseteq \mathbf{G}$ such that $\mathcal{O}_f = \{g^{(\mathbf{u})}{\cdot}f \}_{\mathbf{u} \in S}$. Thus if $\mathbf{u}, \mathbf{v} \in S$ and $g^{(\mathbf{u})}{\cdot}f = g^{(\mathbf{v})}{\cdot}f$ then $\mathbf{u} = \mathbf{v}$.

First of all recall that $e_{\mathbf{0}} = \Lambda = (1/n)(\sum_{\mathbf{m} \in \mathbf{G}} g^{(\mathbf{m})})$ is the integral for $H = k[G]$ determined by $\epsilon(\Lambda) = 1$. For $\mathbf{u} \in S$ there are $|G|/|\mathcal{O}_f| = n/r$ elements $\mathbf{x} \in \mathbf{G}$ such that $g^{(\mathbf{x})}{\cdot}f = g^{(\mathbf{u})}{\cdot}f$. Therefore $(g^{(\mathbf{u})}{\cdot}f)(\Lambda{\cdot}f) = (1/r)(g^{(\mathbf{u})}{\cdot}f)$

For $\mathbf{u}, \mathbf{v} \in S$ we set
$$
E_{\mathbf{u} \; \mathbf{v}} = r(g^{(\mathbf{u})}{\cdot}f \# e_{\mathbf{m}})(g^{(\mathbf{v})}{\cdot}f \# 1).
$$
By definition $E_{\mathbf{u} \; \mathbf{v}} \in \mathcal{M}_{f, \mathbf{m}}$. For $\mathbf{u}, \mathbf{v}, \mathbf{r}, \mathbf{s} \in S$ the calculation
\begin{eqnarray*}
E_{\mathbf{u} \; \mathbf{v}}E_{\mathbf{r} \; \mathbf{s}} & = & r^2(g^{(\mathbf{u})}{\cdot}f \# e_{\mathbf{m}})(g^{(\mathbf{v})}{\cdot}f \# 1)(g^{(\mathbf{r})}{\cdot}f \# e_{\mathbf{m}})(g^{(\mathbf{s})}{\cdot}f \# 1) \\
& = & r^2(g^{(\mathbf{u})}{\cdot}f \# e_{\mathbf{m}})((g^{(\mathbf{v})}{\cdot}f)(g^{(\mathbf{r})}{\cdot}f) \# e_{\mathbf{m}})(g^{(\mathbf{s})}{\cdot}f \# 1) \\
& = & \delta_{\mathbf{v}, \mathbf{r}}r^2(g^{(\mathbf{u})}{\cdot}f \# e_{\mathbf{m}})(g^{(\mathbf{v})}{\cdot}f \# e_{\mathbf{m}})(g^{(\mathbf{s})}{\cdot}f \# 1) \\
& = & \delta_{\mathbf{v}, \mathbf{r}}r^2((g^{(\mathbf{u})}{\cdot}f)((\Lambda g^{(\mathbf{v})}){\cdot}f) \# e_{\mathbf{m}})(g^{(\mathbf{s})}{\cdot}f \# 1) \\
& = & \delta_{\mathbf{v}, \mathbf{r}}r^2((g^{(\mathbf{u})}{\cdot}f)(\Lambda{\cdot}f) \# e_{\mathbf{m}})(g^{(\mathbf{s})}{\cdot}f \# 1) \\
& = & \delta_{\mathbf{v}, \mathbf{r}}r((g^{(\mathbf{u})}{\cdot}f) \# e_{\mathbf{m}})(g^{(\mathbf{s})}{\cdot}f \# 1) \\
& = &  \delta_{\mathbf{v}, \mathbf{r}}E_{\mathbf{u} \; \mathbf{s}}
\end{eqnarray*}
shows that $E_{\mathbf{u} \; \mathbf{v}}E_{\mathbf{r} \; \mathbf{s}} = \delta_{\mathbf{v}, \mathbf{r}}E_{\mathbf{u} \; \mathbf{s}}$.

Our linear isomorphism will be established once we show that the $E_{\mathbf{u} \; \mathbf{v}}$'s form a linearly independent set. Suppose that $\sum_{\mathbf{u}, \mathbf{v} \in S}x_{\mathbf{u}, \mathbf{v}}E_{\mathbf{u} \; \mathbf{v}} = 0$, where $x_{\mathbf{u}, \mathbf{v}} \in k$. Let $\mathbf{r}, \mathbf{s} \in S$. Multiplying both sides of the preceding equation by $g^{(\mathbf{r})}{\cdot}f \# 1$  on the left and $g^{(\mathbf{s})}{\cdot}f \# 1$ on the right gives $x_{\mathbf{r}, \mathbf{s}}E_{\mathbf{r} \; \mathbf{s}} = 0$. To show $x_{\mathbf{r}, \mathbf{s}} = 0$ we need only show $E_{\mathbf{r} \; \mathbf{s}} \neq 0$. This is the case since $E_{\mathbf{r} \; \mathbf{s}}(1 \# e_{\mathbf{m}}) = g^{(\mathbf{r})}{\cdot}f \# e_{\mathbf{m}} \neq 0$.

We leave the justification of the formula for $E_{\mathbf{u} \; \mathbf{v}} \in \mathcal{M}_{f, \mathbf{m}}$ given by
\begin{equation}\label{EqEuv}
E_{\mathbf{u} \; \mathbf{v}}  = \sum_{\mathbf{t} \in \mathbf{G}}\lambda^{((\mathbf{m} - \mathbf{t})(\mathbf{u} - \mathbf{v}))}\left(\sum_{\mathbf{x} \in \mathbf{N}_f}\lambda^{((\mathbf{m} - \mathbf{t})\mathbf{x})}r/n\right) g^{(\mathbf{u})}{\cdot}f \# e_{\mathbf{t}},
\end{equation}
where
$$
\mathbf{N}_f = \{\mathbf{x} \in \mathbf{G} \, | \, g^{(\mathbf{x})}{\cdot}f = f\},
$$
as an exercise for the reader. Note that $|\mathbf{N}_f|1_k \neq 0$ since $|\mathbf{N}_f|$ divides $|G| = n$ and $k$ has a primitive $n^{th}$ of unity by assumption. We will improve on (\ref{EqEuv}).

For $\mathbf{t} \in \mathbf{G}$ let $S_\mathbf{t} = \sum_{\mathbf{x} \in \mathbf{N}_f}\lambda^{((\mathbf{m} - \mathbf{t})\mathbf{x})}$ and let $\mathbf{y} \in \mathbf{N}_f$. Since $\mathbf{N}_f$ is an additive group
$S_\mathbf{t} = \sum_{\mathbf{x} \in \mathbf{N}_f}\lambda^{((\mathbf{m} - \mathbf{t})(\mathbf{x} + \mathbf{y}))} = \lambda^{((\mathbf{m} - \mathbf{t})\mathbf{y})}S_\mathbf{t}$. Thus $S_\mathbf{t} \neq 0$ if and only if $\lambda^{((\mathbf{m} - \mathbf{t})\mathbf{y})} = 1$ for all $\mathbf{y} \in \mathbf{N}_f$ in which case $S_\mathbf{t} = |\mathbf{N}_f|1_k = (n/r)1_k$. Set
$$
\mathbf{I}_f = \{\mathbf{z} \in \mathbf{G} \, | \, \lambda^{(\mathbf{z}\mathbf{y})} = 1 \;\; \mbox{for all $\mathbf{y} \in \mathbf{N}_f$}.\}
$$
Then $\mathbf{I}_f$ is a subgroup of $\mathbf{G}$. Observe that $S_\mathbf{t} \neq 0$ if and only if $\mathbf{m} - \mathbf{t} \in \mathbf{I}_f$. We have shown
\begin{equation}\label{EqEuvBEST}
E_{\mathbf{u} \; \mathbf{v}}  = \sum_{\mathbf{z} \in \mathbf{I}_f}\lambda^{(\mathbf{z}(\mathbf{u} - \mathbf{v}))}g^{(\mathbf{u})}{\cdot}f \# e_{\mathbf{m} - \mathbf{z}}
\end{equation}
for $E_{\mathbf{u} \; \mathbf{v}} \in \mathcal{M}_{f, \mathbf{m}}$.

Let $\mathbf{z}_0 \in \mathbf{I}_f$. Then $E_{\mathbf{u} \; \mathbf{v}}(1 \# e_{\mathbf{m} - \mathbf{z}_0}) = \lambda^{(\mathbf{z}_0(\mathbf{u} - \mathbf{v}))}g^{(\mathbf{u})}{\cdot}f \# e_{\mathbf{m} - \mathbf{z}_0}$. Since $\mathcal{M}_{f, \mathbf{m}}$ is a right ideal of $B \# k[G]$, the linearly independent set
$$
\{g^{(\mathbf{u})}{\cdot}f \# e_{\mathbf{m} - \mathbf{z}} \}_{\mathbf{u} \in S, \, \mathbf{z} \in \mathbf{I}_f}
$$
lies in $\mathcal{M}_{f, \mathbf{m}}$. On the other hand $\mathcal{M}_{f, \mathbf{m}}$ lies in the span of this set by (\ref{EqEuvBEST}). We have shown this linearly independent set is a basis for $\mathcal{M}_{f, \mathbf{m}}$; in particular $|\mathbf{I}_f| = r$.

Observe that $\mathbf{N}_f = \mathbf{N}_{g^{(\mathbf{w})}{\cdot}f}$, and thus  $\mathbf{I}_f = \mathbf{I}_{g^{(\mathbf{w})}{\cdot}f}$, for all $\mathbf{w} \in \mathbf{G}$. This means $\mathbf{N}_f$ and $\mathbf{I}_f$ depend only on the orbit $\mathcal{O}_f$ of $f$ and not a particular element of it.

We next find finding necessary and sufficient conditions for two of these ideals to be the same. Suppose that $\mathcal{M}_{f, \mathbf{m}} = \mathcal{M}_{f', \mathbf{m}'}$, let $L = \mathrm{sp}(\mathcal{O}_f)$, and let $L' = \mathrm{sp}(\mathcal{O}_{f'})$. Either $L = L'$ or $L'L = (0)$. As $L^2 = L$, by (\ref{EqMfm}) we have
$$
\mathcal{M}_{f, \mathbf{m}} = \mathcal{M}_{f, \mathbf{m}}(L \# 1) = \mathcal{M}_{f', \mathbf{m}'}(L \# 1) = (L' \# e_{\mathbf{m}'})(L'L \# 1)
$$
and consequently $L'L \neq (0)$. Therefore $L = L'$. Now it is easy to see that $\mathbf{m} + \mathbf{I}_f = \mathbf{m}' + \mathbf{I}_f$. The two necessary conditions are clearly sufficient.
\begin{Prop}\label{MinimalChar}
Assume the notation above. Then $\mathcal{M}_{f, \mathbf{m}} = \mathcal{M}_{f', \mathbf{m}'}$ if and only if $\mathcal{O}_f = \mathcal{O}_{f'}$ and $\mathbf{m} + \mathbf{I}_f = \mathbf{m}' + \mathbf{I}_f$. \qed
\end{Prop}
\begin{Theorem}\label{ThmAlgebraStructure}
Let $G$ be an abelian group of order $n$,  suppose that the field $k$ contains a primitive $n^{th}$ root of unity, and let $B$ be a left $H$-module algebra which is a power of $k$. Let $H = k[G]$ and let $\mathcal{F}$ be an orthogonal basis of idempotents for $B$. Then:
\begin{enumerate}
\item[{\rm (a)}] $g{\cdot}\mathcal{F} = \mathcal{F}$ for all $g \in G$.
\item[{\rm (b)}] $B \# H$ is a semisimple algebra. Specifically, let $t$ be the number of $G$-orbits of $\mathcal{F}$ and $n_1, \ldots, n_t$ be their respective lengths. Then $n_\ell$ divides $|G|$ for all $1 \leq \ell \leq t$ and
$$B \# H \simeq \bigoplus_{\ell = 1}^t \left(\frac{|G|}{n_\ell}\right){\cdot}\mathrm{M}_{n_\ell}(k).$$
\end{enumerate}
\end{Theorem}

\pf
We need only add a few more remarks to our discussion of minimal ideals of $B \# H$. Let $f \in \mathcal{F}$ and $L = \mathrm{sp}(\mathcal{O}_f)$. Then Since $L$ is an ideal and left $H$-submodule of $B \# H$ it follows that $\sum_{\mathbf{m} \in \mathbf{G}}L \# e_\mathbf{m} = L \# H$ is an ideal $B \# H$. Therefore $L \# H = \sum_{\mathbf{m} \in \mathbf{G}} \mathcal{M}_{f, \mathbf{m}}$ and is the direct sum of $m$ of these ideals. Let $r = |\mathcal{O}_f|$. Then $mr^2 = \mathrm{Dim}(L \# H) = r|G|$. Thus $r$ divides $|G|$ and $m = |G|/r$.

Let $\mathcal{O}_{f_1}, \ldots, \mathcal{O}_{f_t}$ be a list of the distinct $G$-orbits of $\mathcal{F}$ and set $L_i = \mathrm{sp}(\mathcal{O}_{f_i})$ for each $1 \leq i \leq t$. Then $B \# H = \bigoplus_{i = 1}^t (L_i \# H)$ since $B =  \bigoplus_{i = 1}^t L_i$.
\qed
\medskip

Apropos of the preceding theorem, the semisimplicity of $B \# k[G]$ follows from a result of \cite{FM1978}, and from a more general result found in \cite{CF1986}.
\section{The algebra $\mathrm{A}(\mathcal{G}, \theta)$ when $k[\mathcal{G}]$ is a power of $k$.}
Here we consider the algebra structure of $\mathrm{A}(\mathcal{G}, \theta)$ given in (\ref{EqAGTheta}), namely $\mathrm{A}(\mathcal{G}, \theta) \simeq (k[\mathcal{G}] \otimes k[U]) \# k[(\theta)]$, where $k[\mathcal{G}]$ is a power of $k$. Recall that $k[\mathcal{G}]$ is a left $k[(\theta)]$-module algebra, where the module action of $\theta$ on $k[\mathcal{G}]$ is determined by evaluation on $\mathcal{G}$, and that $k[\mathcal{G}] \otimes k[U]$ is a left $k[(\theta)]$-module algebra, where the action is on the first tensorand.

Let $n$ be the order of $\theta$. We will assume that $k$ has a primitive $n^{th}$ root of unity $\lambda$. Since $U$ is identified with the subgroup of $k^\times$ generated by the eigenvalues of $\theta_{k[\mathcal{G}]}$ it follows that $U$ is cyclic and its order divides $n$. Therefore $k$ contains a primitive $|U|^{th}$ root of unity. In particular $k[(\theta)]$ and $k[U]$ are powers of $k$ also.

The algebra $\mathrm{A}(\mathcal{G}, \theta) \simeq (k[\mathcal{G}] \otimes k[U]) \# k[(\theta)]$ is an example of a more general smash product $B \# H$, where $B = C \otimes C'$ is the tensor product of powers of $k$, $H = k[G]$, where $G = (g)$ is a cyclic group of order $n$, $C$ is a left $H$-module algebra, and $B$ is a left $H$-module algebra where $h{\cdot}(c \oplus c') = h{\cdot}c \otimes c'$ for all $h \in H$, $c \in C$, and $c' \in C'$. We will describe the algebra structure of $B \# H$ explicitly.

Let $\mathcal{F}$ and $\mathcal{F}'$ be linear bases for $C$ and $C'$ respectively which are orthogonal families of idempotents. Then $\{f \otimes f'\}_{f \in \mathcal{F}, f' \in \mathcal{F}'}$ is a basis for $C \otimes C'$ of the same type. Thus $C \otimes C'$ is a semisimple algebra whose minimal ideals are the one-dimensional subspaces spanned these basis elements. Let $\{e_i\}_{i \in \mathbb{Z}_n}$ be the basis for $k[G]$ described by (\ref{EqeiOrthogonal}). By results of the preceding section the minimal ideals of $B \# H$ are of the form $\mathcal{M}_{{f \otimes f'}, m}$, where $m \in \mathbb{Z}_n$, and the smash product is the direct sum of the different ones.

Now $\mathcal{M}_{{f \otimes f'}, m} = (L \# e_m)(L \# 1)$, where $L = \mathrm{sp}(\mathcal{O}_{f \otimes f'})$. Let $r = |\mathcal{O}_f|$. Then $\mathcal{O}_f = \{f, g{\cdot}f, \ldots g^{r-1}{\cdot}f\}$ and $g^r{\cdot}f = f$. Thus $r = |\mathcal{O}_{f \otimes f'}|$ and $\mathcal{O}_{f \otimes f'} = \{f \otimes f', g{\cdot}f \otimes f', \ldots g^{r-1}{\cdot}f \otimes f'\}$. For $u, v \in \mathbb{Z}_n$ by (\ref{EqEuvBEST}) we have
$$
E_{u \, v} = \sum_{\ell = 0}^{r-1} \lambda^{(n/r)\ell (u - v)} (g^u{\cdot}f \otimes f') \# e_{m - (n/r)\ell}.
$$
Observe $\mathcal{M}_{f \otimes f', m}$ has basis $\{(g^i{\cdot}f \otimes f') \# e_{m - (n/r)j} \}_{0 \leq i, j \leq r-1}$. The subscripts of the $e$'s form a coset of the $r$-element subgroup of $\mathbb{Z}_n$.

Using the results of Sections \ref{SecConstRev} -- \ref{SecLast} the reader can give detailed descriptions of many Hopf algebras described in Examples \ref{Ex1} -- \ref{Ex6}, notably of the non-trivial Hopf algebras of dimension $12$ where the characteristic of $k$ is not $2$ or $3$ and $k$ is algebraically closed.
\section{The bipoduct of Theorem \ref{ThmCoalgebraStructure} and lower sovability, cosolvability}
The biproduct $A = B \times H$ of Theorem \ref{ThmCoalgebraStructure} satisfies the following properties: $H = k[G]$ is a group algebra of a group $G$ over $k$;
$B$ is a Hopf algebra of ${}_H^H\mathcal{YD}$; there is a normal subgroup $U$ of $G$ which acts trivially on $B$; and $\rho(B) \subseteq k[U] \otimes B$.

Assume that $A$ be any biproduct which satisfies these properties. Since $\rho(B) \subseteq k[U] \otimes B$ we may regard $B$ as a Hopf algebra of ${}_{k[U]}^{k[U]}\mathcal{YD}$. Observe that the inclusion map $B \times k[U] \longrightarrow B \times k[G]$ is a one-one map of Hopf algebras. Let $A_1 = B \times k[U]$. Since $k[U]$ is a normal Hopf subalgebra of $H$ and $H$ is cocommutative, $A_1$ is a normal Hopf subalgebra of $A = B \times H$. The reader is left with the exercise of establishing normality; see (\ref{EqAntipode}) in particular.

Now $A_1^+A = (B^+ \times k[U]) + (B \times k[U]^+H)$ which is the kernel of the onto composition of Hopf algebras maps $B \times H \longrightarrow H \longrightarrow H/k[U]^+H$, where the first is given by $b \times h \mapsto \epsilon(b)h$ and the second is the canonical projection which takes an element to its coset. Thus $A/A_1^+A \simeq H/k[U]^+H$ as Hopf algebras. The group homomorphism $G \longrightarrow G/U$ lifts to an onto Hopf algebra map $H \longrightarrow k[G/U]$ with kernel $k[U]^+H$. Therefore $A/A_1^+A \simeq k[G/U]$ as Hopf algebras.

We have shown $k \subseteq A_1 \subseteq A$ is a lower normal series for $A$ and that the factor $A/A_1^+A \simeq k[G/U]$. Since $k[U]$ acts trivially on $B$ it follows that as algebras $A_1 = B \# k[U] = B \otimes k[U]$. Thus $A_1$ is commutative whenever $B$ is and $U$ is abelian. If this is the case and $G/U$ is abelian then $A$ is lower-solvable; see \cite{MW} for the definition and discussion of lower solvability and related concepts. Our argument has shown:
\begin{Lemma}\label{LemmaLast}
For the biproduct of Theorem \ref{ThmCoalgebraStructure} the sequence $k \subseteq A_1 \subseteq A$ is a lower normal series for $A = k[\mathcal{G}] \times k[G]$, where $A_1 = k[\mathcal{G}] \times k[\mathbf{U}]$. The factor $A/A_1^+A \simeq k[G/\mathbf{U}]$ as Hopf algebras. The factor $A_1$ is a commutative if and only if $\mathcal{G}$ and $U$ are abelian.
\qed
\end{Lemma}

We examine $A_1 = k[\mathcal{G}] \times k[U]$, where $U = \mathbf{U}$, in the context of Theorem \ref{ThmCoalgebraStructure}. In this case the following are satisfied: $\mathcal{G}$ is any group; $U$ is an abelian group; $\theta \in \mathrm{Aut}_{\mathrm{Group}}(\mathcal{G})$, $\theta_{k[\mathcal{G}]}$ is diagonalizable, and the $(\theta)$-orbits of $\mathcal{G}$ are finite; and $k[\mathcal{G}]$ is a Hopf algebra of ${}_{k[U]}^{k[U]}\mathcal{YD}$, where $k[U]$ acts trivially on $k[\mathcal{G}]$ and spans of $(\theta)$-orbits of $\mathcal{G}$ are left $k[U]$-subcomodules of $k[\mathcal{G}]$. We proceed under these assumptions.

First of all suppose $\mathcal{N}$ is a normal subgroup of $\mathcal{G}$ and $\theta(\mathcal{N}) \subseteq \mathcal{N}$. The inclusion implies that $\mathcal{N}$ is the union of $(\theta)$-orbits. Therefore $k[\mathcal{N}]$ is a left $k[U]$-subcomodule of $k[\mathcal{G}]$. Consequently $k[\mathcal{N}]$ is a Hopf algebra of ${}_{k[U]}^{k[U]}\mathcal{YD}$. We regard the biproduct $k[\mathcal{N}] \times k[U]$ as a Hopf subalgebra of $k[\mathcal{G}] \times k[U]$ in the usual way; indeed it is a normal Hopf subalgebra since $k[\mathcal{N}]$ is a normal Hopf subalgebra of $k[\mathcal{G}]$. The normality assertions are more easily verified using the fact that $k[\mathcal{G}] \#  k[U] = k[\mathcal{G}] \otimes k[U]$ as algebras.

Let $A_2 = k[\mathcal{N}] \times k[U]$. Then $A^+_2A_1 = k[\mathcal{N}]^+k[\mathcal{G}] \otimes k[U] + k[\mathcal{G}] \otimes k[U]^+$, the kernel of the onto composite $k[\mathcal{G}] \times k[U] \longrightarrow k[\mathcal{G}] \longrightarrow k[\mathcal{G}]/k[\mathcal{N}]^+k[\mathcal{G}]$ of Hopf algebra maps, where the first is determined by $b \times h \mapsto \epsilon(h)b$ and the second is the canonical projection.  Thus $A_1/A_2^+A_1 \simeq k[\mathcal{G}]/k[\mathcal{N}]^+k[\mathcal{G}] \simeq k[\mathcal{G}/\mathcal{N}]$ as Hopf algebras.

The commutator subgroup $\mathcal{G}^{(1)}$ of $\mathcal{G}$ is a fully invariant normal subgroup of $\mathcal{G}$; in particular $\theta(\mathcal{G}^{(1)}) = \mathcal{G}^{(1)}$. With $\mathcal{N} = \mathcal{G}^{(1)}$ we have $A_2 = k[\mathcal{G}^{(1)}] \times k[U]$ is a normal Hopf subalgebra of $A_1$ and $A_1/A_2^+A_1 \simeq k[\mathcal{G}/\mathcal{G}^{(1)}]$. By induction we can construct as sequence of Hopf subalgebras
$$
A_2 \supseteq A_3 \supseteq A_4 \supseteq \cdots
$$
where $A_\ell = k[\mathcal{G}^{(\ell - 1)}] \times k[U]$, $A_{\ell +1}$ is a normal Hopf subalgebra of $A_\ell$, and $A_\ell/A_{\ell +1}^+A_\ell \simeq k[\mathcal{G}^{(\ell - 1)}/\mathcal{G}^{(\ell)}]$ for all $\ell \geq 1$. Here $\mathcal{G}^{(0)} = \mathcal{G}$ and $\mathcal{G}^{(\ell + 1)}$ is the commutator subgroup of $\mathcal{G}^{(\ell)}$ for $\ell \geq 0$. In light of Lemma \ref{LemmaLast} the first statement of the following is established.

\begin{Theorem}\label{ThmLast}
The biproduct of Theorem \ref{ThmCoalgebraStructure} is lower solvable when $G$ is abelian and $\mathcal{G}$ is solvable. Thus the biproduct is lower solvable when $G$ is abelian and $\mathcal{G}$ is finite and has odd order.
\end{Theorem}

\pf
We need only note that finite groups of odd order are solvable by the Feit-Thompson Theorem \cite[Chapter 1, \S 1]{FT}.
\qed
\begin{Theorem}\label{ThmLastLast}
The biproduct of Theorem \ref{ThmCoalgebraStructure} is lower cosolvable.
\end{Theorem}

\pf
Let $A = B \times k[G]$ be the biproduct of Theorem \ref{ThmCoalgebraStructure}, where $B = k[\mathcal{G}]$. In light of Lemma \ref{LemmaLast} we need only show that $A_1 = B \times k[\mathbf{U}]$ is cosolvable. Since $\mathcal{H} = k1 \times k[\mathbf{U}]$ is in the center of $A_1$ the adjoint actions of $A_1$ on $\mathcal{H}$ are trivial. Thus $\mathcal{H}$ is a normal Hopf subalgebra of $A_1$.

We next show that $A_1/\mathcal{H}^+A_1 \simeq k[\mathcal{G}]$ as Hopf algebras. It is easy to see $\mathcal{H}^+A_1 = k[\mathcal{G}] \times k[\mathbf{U}]^+$. The algebra map $f :  k[\mathcal{G}] \otimes k[\mathbf{U}] \longrightarrow k[\mathcal{G}]$, where $f =  \mathrm{I}_{k[\mathcal{G}]}  \otimes \epsilon$, has kernel $\mathcal{H}^+A_1$. We need only show that $f$ is a coalgebra map; for then we can conclude the induced map $A_1/\mathcal{H}^+A_1 \simeq k[\mathcal{G}]$ is an isomorphism of Hopf algebras.

 Now $A_1$ is the sum of subcoalgebras of the form $\mathrm{C}(\mathcal{O}_b, \mbox{\boldmath $\lambda$})(1 \times g)$, where $b \in \mathcal{G}$ and $g \in G$. The subcoalgebra $\mathrm{C}(\mathcal{O}_b, \mbox{\boldmath $\lambda$})$ has a basis $\{c_{i \, j}\}_{0 \leq i \leq r-1}$, where $r = |\mathcal{O}_b|$, $c_{i \, j} = \sum_{\ell = 0}^{r-1} \theta^i(b) \times (\lambda^{\ell (i - j)}/r) \mbox{\boldmath{$\lambda$}${}^\ell$}$, and $\lambda \in k^\times$ is a primitive $r^{th}$ root of unity. Thus $f(c_{i \, j}) = \delta_{i, j}\theta^i(b)$ for all $0 \leq i, j \leq r-1$. Note that $f(1 \times g) = 1$.

 Our calculations show that the algebra maps $(f \otimes f){\circ}\Delta$, $\Delta{\circ} f$ and $\epsilon{\circ}f$, $\epsilon$ agree on generators for $A_1$. Therefore $(f \otimes f){\circ}\Delta = \Delta{\circ} f$ and $\epsilon{\circ}f = \epsilon$ which means $f$ is a coalgebra map.

In light of Lemma \ref{LemmaLast} we have shown $k \subseteq \mathcal{H} \subseteq A_1 \subseteq A$ is a lower normal series for $A$ whose factors are group algebras, thus are cocommutative. By definition $A$ is lower cosolvable.
\qed
\section{Non-trivial semisimple cosemisimple Hopf algebras with a unique proper normal Hopf subalgebra}
Throughout this section $A = k[\mathcal{G}] \times k[G]$ is the biproduct of Theorem \ref{ThmCoalgebraStructure}, where $G = \mathrm{Ker}(\pi)$ and is abelian. Set $B = k[\mathcal{G}]$ and $H = k[G]$. We will make heavy use of the discussion of Section \ref{SecConstRev} and freely use the results therein. Since the $H$-module action on $B$ is trivial $B \times H = B \otimes H$ as algebras. Since $H$ is also commutative, the adjoint actions of $A$ on itself are given by
\begin{equation}\label{EqAdjointALeft}
(b \times h){\succ} (b' \times h') = (b{\succ}b') \times \epsilon(h)h'
\end{equation}
and
\begin{equation}\label{EqAdjointARight}
 (b' \times h'){\prec} (b \times h) = (b'{\prec}b_{(0)}) \times \epsilon(h)S(b_{(-1)})h'
\end{equation}
for all $b, b' \in B$ and $h, h' \in H$. For these calculations the formula
$$
S(b \times h) = S(b_{(0)}) \times S(b_{(-1)}h)
$$
for all $b \in B$ and $h \in H$ is handy. We use $S$ to denote the antipode of $A$, $B$, or $H$; which one should always be clear from context.

Let $\mathcal{A}$ be a Hopf subalgebra of $A$. We first show that $\mathcal{A}$ has a basis consisting of elements of the form $b \times g$, where $b \in \mathcal{G}$ and $g \in G$. Properties of these pairs $(b, g)$ will guide of analysis of $\mathcal{A}$. Observe that the set of tensors $\mathcal{G} \times G$ is the set $\mathrm{G}(B \otimes H)$ when $B \otimes H$ is viewed as the Hopf algebra which is the tensor product of group algebras over $k$.

First a notation simplification. We write $\mathrm{C}(\mathcal{O})$ for $\mathrm{C}(\mathcal{O}, \mathbf{U}_{|\mathcal{O}|})$ since the subgroup $\mathbf{U}_{|\mathcal{O}|}$ of $\mathbf{U}$ is determined by $r = |\mathcal{O}|$.

Since $A$ is cosemisimple $\mathcal{A}$ is as well. Therefore $\mathcal{A}$ is the direct sum of simple subcoalgebras of the form $\mathrm{C}(\mathcal{O})(1 \times g)$, where $\mathcal{O} = \mathcal{O}_b$ for some $b \in \mathcal{G}$ and $g \in G$. The set of tensors $\mathcal{O}_b \times \mathbf{U}_rg$, where $r = |\mathcal{O}_b|$, is a basis for $\mathrm{C}(\mathcal{O})(1 \times g)$. Therefore $\mathcal{A}$ has a basis consisting of elements of the form $b \times g$, where $b \in \mathcal{G}$ and $g \in G$. We denote this basis by $\mathcal{B}_\mathcal{A}$.

Suppose $b \times g \in \mathcal{B}_\mathcal{A}$. Then $b \times g \in \mathcal{A} \cap (\mathrm{C}(\mathcal{O}_b)(1 \times g))$ which means $\mathrm{C}(\mathcal{O}_b)(1 \times g) \subseteq \mathcal{A}$ since the former is a simple subcoalgebra of $A$. Since $\mathcal{A}$ is a Hopf subalgebra of $A$ it follows that $S(\mathrm{C}(\mathcal{O}_b)(1 \times g)) \subseteq \mathcal{A}$. Now
$$
S(\mathrm{C}(\mathcal{O}_b)(1 \times g)) = (1 \times g^{-1})S(\mathrm{C}(\mathcal{O}_b)) = S(\mathrm{C}(\mathcal{O}_b))(1 \times g^{-1})
$$
since $H$ is commutative. Also $S(\mathcal{O}_b) = \mathcal{O}_{S(b)} = \mathcal{O}_{b^{-1}}$ since $S{\circ}\theta = \theta{\circ}S$. In particular $\mathbf{U}_{|\mathcal{O}_{b^{-1}}|} = \mathbf{U}_{|\mathcal{O}_b|}$. By assumption there is a primitive $r^{th}$ root of unity $\lambda \in k^\times$. Now $\mathbf{U}_{|\mathcal{O}_b|} = (\mbox{\boldmath $\lambda$})$ and $\{b_{\lambda^\ell}\}_{0 \leq i \leq r-1}$ is a basis for $\mathrm{sp}(\mathcal{O}_b)$. Likewise $\{b^{-1}_{\lambda^\ell}\}_{0 \leq i \leq r-1}$ is a basis for $\mathrm{sp}(\mathcal{O}_{b^{-1}})$. Since $S(b_{\lambda^\ell}) = b^{-1}_{\lambda^\ell}$, for $u \in  \mathbf{U}_{|\mathcal{O}_b|}$ we have
$$
S(b_{\lambda^\ell} \times u) = S(b_{\lambda^\ell}) \times S(\mbox{\boldmath $\lambda$}^\ell u) \in \mathrm{sp}(\mathcal{O}_{b^{-1}}) \times \mathbf{U}_{|\mathcal{O}_{b^{-1}}|} = \mathrm{C}(\mathcal{O}_{b^{-1}}).
$$
We have shown that $S(\mathrm{C}(\mathcal{O}_b)) \subseteq \mathrm{C}(\mathcal{O}_b^{-1})$; thus the preceding subspaces are equal and as a result $S(\mathrm{C}(\mathcal{O}_b)(1 \times g)) = \mathrm{C}(\mathcal{O}_{b^{-1}})(1 \times g^{-1})$. In particular $b^{-1} \times g^{-1} \in \mathcal{B}_\mathcal{A}$.
\begin{Lemma}\label{LemmaMathcalA}
Let $A = k[\mathcal{G}] \times k[G]$ be a biproduct described in Theorem \ref{ThmCoalgebraStructure}, where $G = \mathrm{Ker}(\pi)$ and is abelian. Suppose $\mathcal{A}$ is a Hopf subalgebra of $A$. Then:
\begin{enumerate}
\item[{\rm (a)}] $1 \times 1 \in \mathcal{B}_\mathcal{A}$.
\end{enumerate}
Suppose $b \times g \in \mathcal{B}_\mathcal{A}$. Then:
\begin{enumerate}
\item[{\rm (b)}] $\theta^i(b) \times ug \in \mathcal{B}_\mathcal{A}$ for all $i \geq 0$ and $u \in \mathbf{U}_r$, where $r = |\mathcal{O}_b|$.
\item[{\rm (c)}] $b^{-1} \times g^{-1} \in \mathcal{B}_\mathcal{A}$.
\item[{\rm (d)}] $1 \times \mathbf{U}_r \subseteq \mathcal{B}_\mathcal{A}$, where $r = |\mathcal{O}_b|$.
\item[{\rm (e)}] If $b' \times g' \in \mathcal{B}_\mathcal{A}$ also then $bb' \times gg' \in \mathcal{B}_\mathcal{A}$.
\end{enumerate}
\qed
\end{Lemma}

\pf
Parts (b) and (c) result from the discussion preceding the lemma. Parts (a) and (e) follow since $\mathcal{A}$ is a subalgebra of $A$. Part (d) follows from parts (b), (c), and (e).
\qed
\medskip

We next consider adjoint actions on and by $\mathcal{A}$.
\begin{Lemma}\label{LemmaAdjointA}
Let $A = k[\mathcal{G}] \times k[G]$ be a biproduct described in Theorem \ref{ThmCoalgebraStructure}, where $G = \mathrm{Ker}(\pi)$ and is abelian. Suppose $\mathcal{A}, \mathcal{A}'$ are Hopf subalgebras of $A$ and $\mathcal{A}'  \subseteq \mathcal{A}$. Then:
\begin{enumerate}
\item[{\rm (a)}] $\mathcal{A} {\succ} \mathcal{A}' \subseteq \mathcal{A}'$ if and only if $(bb'b^{-1}) \times g' \in \mathcal{B}_{\mathcal{A}'}$ for all $b \times g \in \mathcal{B}_\mathcal{A}$ and $b' \times g' \in \mathcal{B}_{\mathcal{A}'}$.
\item[{\rm (b)}] $\mathcal{A}'  {\prec} \mathcal{A} \subseteq \mathcal{A}'$ if and only if $(b' {\prec} b_{\lambda^\ell}) \times \mbox{\boldmath $\lambda$}^{-\ell}g' \in \mathcal{A}'$ for all $b' \times g' \in \mathcal{B}_{\mathcal{A}'}$ and $b \times g \in \mathcal{B}_\mathcal{A}$.
\end{enumerate}
\end{Lemma}

\pf
Part (a) follows directly by (\ref{EqAdjointALeft}). As for part (b), let $b' \times g' \in \mathcal{B}_{\mathcal{A}'}$ and $b \times g \in \mathcal{B}_\mathcal{A}$. By part (b) of Lemma \ref{LemmaMathcalA} we see that $\theta^i(b) \times g \in \mathcal{B}_\mathcal{A}$ for all $i \geq 0$. Let $r = |\mathcal{O}_b|$ and $\lambda \in k^\times$ be a primitive $r^{th}$ root of unity.

Suppose that $\mathcal{A}'  {\prec} \mathcal{A} \subseteq \mathcal{A}'$ and fix $0 \leq i \leq r-1$. Using (\ref{EqAdjointARight}) we see
$$
(b' \times g'){\prec}(\theta^i(b) \times g) = b'{\prec}(\theta^i(b)_{(0)}) \times S(\theta^i(b)_{(-1)})g' \in \mathcal{A}';
$$
thus by (\ref{EqThetaIb}) we have $\sum_{\ell \in \mathbb{Z}_r} (b'{\prec}b_{\lambda^\ell}) \times \lambda^{i\ell}\mbox{\boldmath $\lambda$}^{-\ell}g' \in \mathcal{A}'$, where $(\mbox{\boldmath $\lambda$}) = \mathbf{U}_r$. Let $\mathbf{a}_\ell = (b'{\prec}b_{\lambda^\ell}) \times \mbox{\boldmath $\lambda$}^{-\ell}g'$ for $\ell \in \mathbb{Z}_r$. Denote the $k$-vector space of $1 \times r$ matrices with entries in $\mathcal{A}'$ by $\mathcal{A}'^r$ and let $\mathbf{a} \in \mathcal{A}'^r$ have entries $\mathbf{a}_0, \ldots, \mathbf{a}_{r-1}$. Since $\mathbf{V} = (v_{i \, j}) \in \mathrm{M}_r(\mathcal{A}')$ is invertible, where $v_{i \, j} = (\lambda^{ij}/r)1$, and $\mathbf{a}' := \mathbf{V}\mathbf{a} \in \mathcal{A}'^r$, it follows that $\mathbf{a} = \mathbf{V}^{-1}\mathbf{a}' \in \mathcal{A}'^r$. Therefore $\mathbf{a}_0, \ldots, \mathbf{a}_{r-1} \in \mathcal{A}'$.

We have shown that the condition of part (b) is necessary. It is easy to see that it is also sufficient.
\qed
\medskip

Let
$$
\mathcal{G}_\mathcal{A} = \{ b \in \mathcal{G} \, | \, b \times g \in \mathcal{B}_\mathcal{A} \; \mbox{for some $g \in G$} \}
$$
and
$$
G_\mathcal{A} = \{ g \in G \, | \, b \times g \in \mathcal{B}_\mathcal{A} \; \mbox{for some $b \in \mathcal{G}$} \}.
$$
By Lemma \ref{LemmaMathcalA} we have $\mathcal{G}_\mathcal{A} = \theta(\mathcal{G}_\mathcal{A})$ and is a subgroup of $\mathcal{G}$; also $G_\mathcal{A}$ is a subgroup of $G$. We observe that
\begin{equation}\label{RhoBNA}
\rho(\mathrm{sp}(\mathcal{O}_b)) \subseteq k[N_\mathcal{A}] \otimes \mathrm{sp}(\mathcal{O}_b)
\end{equation}
for all $b \in \mathcal{G}_\mathcal{A}$ which follows by (\ref{RhoB}) and part (d) of Lemma \ref{LemmaMathcalA}. Therefore $k[\mathcal{G}_\mathcal{A}]$ is a left $k[N_\mathcal{A}]$-subcomodule of $k[\mathcal{G}]$. Consequently $k[\mathcal{G}_\mathcal{A}] \times k[G_\mathcal{A}]$ is a Hopf subalgebra of $A$ and $\mathcal{A} \subseteq k[\mathcal{G}_\mathcal{A}] \times k[G_\mathcal{A}]$ by part (a) of Lemma \ref{LemmaMathcalA}.

Let $\iota : H \longrightarrow A$ be the one-one map of Hopf algebras given by $\iota (h) = 1 \times h$ for all $h \in H$. Then $\iota^{-1}(\mathcal{A}) = k[N_\mathcal{A}]$ for some subgroup $N_\mathcal{A}$ of $G$. We emphasize $1 \times N_\mathcal{A} \subseteq \mathcal{B}_\mathcal{A}$.

There is an onto group homomorphism $f_\mathcal{A} : \mathcal{G}_\mathcal{A} \longrightarrow G_\mathcal{A}/N_\mathcal{A}$ defined as follows: $f_\mathcal{A}(b) = gN_\mathcal{A}$, where $b \times g \in \mathcal{B}_\mathcal{A}$. Our assertion follows by Lemma \ref{LemmaMathcalA} once we show that $f_\mathcal{A}$ is well-defined. First of all $N_\mathcal{A}$ is a normal subgroup of $G$ since $G$ is abelian. Let $b \in \mathcal{G}_\mathcal{A}$. Then $b \times g \in \mathcal{B}_\mathcal{A}$ for some $g \in G_\mathcal{A}$. Suppose $b\times g' \in \mathcal{B}_\mathcal{A}$ also. Then $1 \times g^{-1}g' = (b^{-1} \times g^{-1})(b \times g') \in \mathcal{B}_\mathcal{A}$ by parts (c) and (e) of Lemma \ref{LemmaMathcalA}. Therefore $g^{-1}g' \in N_\mathcal{A}$, or equivalently $gN_\mathcal{A} = g'N_\mathcal{A}$. We have shown $f_\mathcal{A}$ is well-defined. Observe that $f_\mathcal{A}{\circ}\theta = f_\mathcal{A}$ by part (b) of Lemma \ref{LemmaMathcalA}.

Let $\mathcal{N}_\mathcal{A} = \mathrm{Ker}(f_\mathcal{A})$. Then $\mathcal{N}_\mathcal{A}$ is a normal subgroup of $\mathcal{G}_\mathcal{A}$. Suppose $b \in \mathcal{N}_\mathcal{A}$. Then $b \times g \in \mathcal{B}_\mathcal{A}$, where $g \in N_\mathcal{A}$. Since $1 \times N_\mathcal{A} \subseteq \mathcal{B}_\mathcal{A}$ we have $b \times N_\mathcal{A} = (b \times g)(1 \times N_\mathcal{A}) \subseteq \mathcal{B}_\mathcal{A}$ by part (e) of Lemma \ref{LemmaMathcalA}. Now $\theta(\mathcal{N}_\mathcal{A}) \subseteq \mathcal{N}_\mathcal{A}$ since $f_\mathcal{A}{\circ}\theta = f_\mathcal{A}$. Therefore $k[\mathcal{N}_\mathcal{A}]$ is a left $k[N_\mathcal{A}]$-subcomodule by (\ref{EqRhoNY}). Consequently $k[\mathcal{N}_\mathcal{A}] \times k[N_\mathcal{A}]$ is a Hopf subalgebra of $A$ and $k[\mathcal{N}_\mathcal{A}] \times k[N_\mathcal{A}] \subseteq \mathcal{A}$. Set
$$
L_\mathcal{A} = k[\mathcal{N}_\mathcal{A}] \times k[N_\mathcal{A}] \;\; \mbox{and} \;\; U_\mathcal{A} = k[\mathcal{G}_\mathcal{A}] \times k[G_\mathcal{A}].
$$
\begin{Prop}\label{PropLastLastLast}
Let $A = k[\mathcal{G}] \times k[G]$ be a biproduct described in Theorem \ref{ThmCoalgebraStructure}, where $G = \mathrm{Ker}(\pi)$ and is abelian. Suppose $\mathcal{A}$ is a Hopf subalgebra of $A$. Then
\begin{enumerate}
\item[{\rm (a)}] $L_\mathcal{A} \subseteq \mathcal{A} \subseteq U_\mathcal{A}$.
\item[{\rm (b)}] $\mathcal{N}_\mathcal{A} = \mathcal{G}_\mathcal{A}$ and only if $N_\mathcal{A} = G_\mathcal{A}$ if and only if $L_\mathcal{A} = U_\mathcal{A}$ if and only if $L_\mathcal{A} = \mathcal{A} = U_\mathcal{A}$.
\item[{\rm (c)}] $A {\succ}\mathcal{A} \subseteq \mathcal{A}$ implies $\mathcal{G}_\mathcal{A}$ is a normal subgroup of $\mathcal{G}$.
\item[{\rm (d)}] $L_\mathcal{A}$ is a normal Hopf subalgebra of $U_\mathcal{A}$ and $\mathcal{A}$, and $\mathcal{A}$ is a normal Hopf subalgebra of $U_\mathcal{A}$.
\item[{\rm (e)}] The resulting quotients $U_\mathcal{A}/L_\mathcal{A}^+U_\mathcal{A}$, $\mathcal{A}/(L_\mathcal{A})^+\mathcal{A}$, and $U_\mathcal{A}/\mathcal{A}^+U_\mathcal{A}$, are commutative Hopf algebras.
\end{enumerate}
\end{Prop}

\pf
Part (a) was noted above; part (b) follows easily by virtue of the group isomorphism $\mathcal{G}_\mathcal{A}/ \mathcal{N}_\mathcal{A} \simeq G_\mathcal{A}/ N_\mathcal{A}$ and part (a). We establish parts (c) and (d) together.

The implicit assertions about the left adjoint actions follow easily from part (a) of Lemma \ref{LemmaAdjointA} and the fact that $\mathcal{N}_\mathcal{A}$ is a normal subgroup of $\mathcal{G}_\mathcal{A}$, with the exception of $U_\mathcal{A} {\succ} \mathcal{A} \subseteq \mathcal{A}$.

Suppose $b \times g \in \mathcal{B}_{U_\mathcal{A}}$ and $b' \times g' \in \mathcal{B}_\mathcal{A}$. We must show $(b \times g){\succ}(b' \times g') = bb'b^{-1} \times g' \in \mathcal{A}$. Since $b \in \mathcal{G}_\mathcal{A}$ by definition $g \times h \in \mathcal{B}_\mathcal{A}$ for some $h \in G_\mathcal{A}$. Therefore $bb'b^{-1} \times g' = (b \times h)(b' \times g')(b^{-1} \times h^{-1}) \in \mathcal{B}_\mathcal{A}$ by parts (c) and (e) of Lemma \ref{LemmaMathcalA}.

The implicit right adjoint action assertions boil down to $L_\mathcal{A} \subseteq L_\mathcal{A} {\prec} U_\mathcal{A}$ and $\mathcal{A} \subseteq \mathcal{A} {\prec} U_\mathcal{A}$.

First we show $L_\mathcal{A} \subseteq L_\mathcal{A} {\prec} U_\mathcal{A}$. Now $\rho(k[\mathcal{G}_\mathcal{A}]) \subseteq k[N_\mathcal{A}] \otimes k[\mathcal{G}_\mathcal{A}]$ by (\ref{RhoBNA}). Since $\mathcal{N}_\mathcal{A}$ is a normal subgroup of $\mathcal{G}_\mathcal{A}$ and $N_\mathcal{A}$ is a subgroup of $G$ we can use (\ref{EqAdjointARight}) to calculate
$$
L_\mathcal{A}{\prec}U_\mathcal{A} \subseteq (k[\mathcal{N}_\mathcal{A}]{\prec} k[\mathcal{G}_\mathcal{A}]) \times k[N_\mathcal{A}]k[N_\mathcal{A}] \subseteq k[\mathcal{N}_\mathcal{A}] \times k[N_\mathcal{A}] = L_\mathcal{A}.
$$
Next we show $\mathcal{A} \subseteq \mathcal{A} {\prec} U_\mathcal{A}$. Let $b' \times g' \in \mathcal{B}_\mathcal{A}$. By (\ref{EqAdjointARight})
$$
(b' \times g'){\prec}U_\mathcal{A} \in (b'{\prec}k[\mathcal{G}_\mathcal{A}]) \times k[N_\mathcal{A}]g' = (k[\mathcal{G}_\mathcal{A}]{\succ}b') \times k[N_\mathcal{A}]g';
$$
thus
$$
(k[\mathcal{G}_\mathcal{A}]{\succ}b') \times k[N_\mathcal{A}]g' \subseteq (U_\mathcal{A}{\succ}(b' \times g'))(1 \times k[N_\mathcal{A}]) \subseteq \mathcal{A}(1 \times k[N_\mathcal{A}]) \subseteq \mathcal{A}.
$$
Therefore $\mathcal{A} \subseteq \mathcal{A} {\prec} U_\mathcal{A}$ and parts (c) and (d) are established.

We show part (e). Let $K = U_\mathcal{A}$ or $K = \mathcal{A}$ and let $b \times g, b' \times g' \in \mathcal{B}_\mathcal{A}$. Then $b, b'  \in \mathcal{G}_\mathcal{A}$ and $g, g' \in G_\mathcal{A}$. Since $f_\mathcal{A} : \mathcal{G}_\mathcal{A} \longrightarrow G_\mathcal{A}/N_\mathcal{A}$ is a group homomorphism to an abelian group, $b'b = nbb'$ for some $n \in \mathrm{Ker}(f_\mathcal{A}) = \mathcal{N}_\mathcal{A}$. Therefore
\begin{eqnarray*}
d & = & (b' \times g')(b \times g) - (b \times g)(b' \otimes g') \\
 & = & ((n - 1) \times 1)(b \times g)(b' \times g') \\
 & \in & L_\mathcal{A}^+ (b \times g)(b' \times g') \\
 & \subseteq  & \mathcal{A}^+(b \times g)(b' \times g').
\end{eqnarray*}
By part (e) of Lemma \ref{LemmaMathcalA} when $K = U_\mathcal{A}$, $d \in L_\mathcal{A}^+U_\mathcal{A}$, $\mathcal{A}^+U_\mathcal{A}$; when $K = \mathcal{A}$, $d \in L_\mathcal{A}^+\mathcal{A}$. We have shown that all three quotients are commutative.
\qed

We come to the main result of this section.
\begin{Theorem}\label{ThmMainOneNormal}
Let $\mathcal{G}$ be a finite non-abelian simple group and $\theta \in \mathrm{Aut}_{\mathrm{Group}}(\mathcal{G})$ have prime order $p$. Assume that $k$ has a primitive $p^{th}$ root of unity. Then there exists a biproduct $A = k[\mathcal{G}] \times k[G]$, where $G$ is cyclic of order $p$, such that:
\begin{enumerate}
\item[{\rm (a)}] $A$ is a non-trivial semisimple cosemisimple Hopf algebra over $k$.
\item[{\rm (b)}] The normal Hopf subalgebras of $A$ are $k1$, $\mathcal{H} = k1 \times k[G]$, and $A$.
\item[{\rm (c)}] $\mathcal{H} \simeq k[G]$ and $A/\mathcal{H}^+A \simeq k[\mathcal{G}]$ as Hopf algebras. In particular $A$ is not lower solvable.
\end{enumerate}
\end{Theorem}

\pf
The endomorphism $\theta_{k[\mathcal{G}]}$ of $k[\mathcal{G}]$ has order $p$. Its eigenvalues belong to the subgroup $U \subseteq k^\times$ of $p^{th}$ roots of unity. Since $U$ and $\theta_{k[\mathcal{G}]}$ have prime order $p$, necessarily $U$ is the set of eigenvectors of $\theta_{k[\mathcal{G}]}$. Thus $\theta_{k[\mathcal{G}]}$ is diagonalizable. We take $G = U$. Our biproduct is that of Proposition \ref{PropLastLastLast}. We have established part (a).

To show part (b), suppose $\mathcal{A}$ is a normal Hopf subalgebra of $A$. Then $\mathcal{G}_\mathcal{A}$ is a normal subgroup of $\mathcal{G}$ by part (c) of Proposition \ref{PropLastLastLast}. Since $\mathcal{G}$ is simple either $\mathcal{G}_\mathcal{A} = \mathcal{G}$ or $\mathcal{G}_\mathcal{A} = \{ 1 \}$. Likewise $G_\mathcal{A} = G$ or $G_\mathcal{A} = \{ 1\}$.

Since $\mathcal{G}_\mathcal{A}/\mathcal{N}_\mathcal{A} \simeq G_\mathcal{A}/N_\mathcal{A}$ and the latter is abelian, $\mathcal{G}_\mathcal{A} = \mathcal{N}_\mathcal{A}$ and therefore $G_\mathcal{A} = N_\mathcal{A}$. Thus $\mathcal{A} = U_\mathcal{A} = k[\mathcal{G}_\mathcal{A}] \times k[G_\mathcal{A}]$ by part (b) of Proposition \ref{PropLastLastLast}.

 Suppose $G_\mathcal{A} = \{ 1\}$. Then $\mathcal{B}_\mathcal{A} = \mathcal{G}_\mathcal{A} \times 1$. Let $b \in \mathcal{G}_\mathcal{A}$. Then $1 \times \mathbf{U}_{|\mathcal{O}_b|} \subseteq \mathcal{B}_\mathcal{A}$ by part (d) of Lemma \ref{LemmaMathcalA}. Hence $\mathbf{U}_{|\mathcal{O}_b|} = \{ 1 \}$ which means $|\mathcal{O}_b| = 1$, or equivalently $\theta(b) = b$. Since $\theta \neq \mathrm{I}_\mathcal{G}$ necessarily $\mathcal{G}_\mathcal{A} \neq \mathcal{G}$. Thus $\mathcal{G}_\mathcal{A} = \{ 1 \}$ from which $\mathcal{A} = k(1 \times 1)$ follows.

 Now suppose $G_\mathcal{A} = G$. Then $\mathcal{A} = A$ or $\mathcal{A} = k1 \times k[G]$. At this point we recall that the normality of the latter and part (c) was established in the proof of Theorem \ref{ThmLastLast}.
 \qed
 \medskip

Any finite non-abelian simple group $\mathcal{G}$ accounts for a biproduct which satisfies the conclusion of the preceding theorem, provided $k$ has enough roots of unity. For $b \in \mathcal{G}$ consider $\theta_b \in \mathrm{Aut}_{\mathrm{Group}}(\mathcal{G})$ defined by $\theta_b(x) = bxb^{-1}$ for all $x \in \mathcal{G}$. Since $\mathcal{G}$ is not abelian $\theta_b \neq \mathrm{I}_\mathcal{G}$ for some $b \in \mathcal{G}$. Choose such a $b$. Then $|(\theta_b)| = n > 1$. For any prime divisor $p$ of $n$ there is a $\theta \in (\theta_b)$ of order $p$. The hypothesis of the preceding theorem is met provided $k$ has a primitive $p^{th}$ root of unity for some such prime $p$.

\end{document}